\documentclass[12pt]{article}
\usepackage{amssymb}
\usepackage{amsthm}\usepackage{epsf,epsfig,amsfonts}
\usepackage{amsmath}
\usepackage{graphicx}
\usepackage{a4wide}
\newtheorem{corollary*}{Corollary}

\def\trace{\mbox{\rm trace}}
\newcommand{\be}{\begin{equation}}
\newcommand{\ee}{\end{equation}}

\newcommand{\weg}[1]{}
\newcommand{\s}[1]{\stackrel{#1}{\tau}}
\newcommand{\C}[1]{\stackrel{#1}{c}}
\usepackage{epsf,epsfig,amsfonts,a4wide}

\usepackage{amsfonts}
\usepackage{amsmath,amssymb,amsthm}
\usepackage{url}
\usepackage{amsmath,amssymb,amsthm}
\usepackage{amsfonts}
\usepackage{amsmath,amssymb,amsthm}
\usepackage{url}

\newtheorem{Th}{Theorem}

\newtheorem{Lemma}{Lemma}
\newtheorem{Cor}{Corollary}

\theoremstyle{remark}

\newtheorem{Rem}{Remark}

\newcommand{\const}{\mbox{\rm const}}
\title{ Proof  of  the projective 
Lichnerowicz conjecture for  pseudo-Riemannian metrics  with degree of mobility greater than  two}
\date{}
\author{ Volodymyr  Kiosak\footnote{Partially supported by DFG (SPP 1154)  }, Vladimir  S. Matveev\footnote{Partially supported by DFG (SPP 1154 and GK 1523)  } \thanks{ Institute of Mathematics, FSU Jena, 07737 Jena Germany,  vladimir.matveev@uni-jena.de} }
\begin{document}
\maketitle

\section{Introduction}
\subsection{Definitions and result} \label{results} 

Let $M$ be a connected manifold of dimension $n\ge 3$, let  $g$ be  a (Riemannian or pseudo-Riemannian) metric on it.  
We say that a metric $\bar g$ on the same manifold $M$ is \emph{geodesically equivalent}  to $g$, if every $g$-geodesic is  a reparametrized 
$\bar g$-geodesic.   We say that they are   \emph{affine  equivalent}, if their Levi-Civita connections coincide.

As we recall in Section \ref{standard},  
the set of metrics geodesically equivalent  to a given one (say, $g$) is  in one-to-one correspondence with the nondegenerate solutions of the equation \eqref{basic}. Since the equation \eqref{basic} is linear, the space of its solutions is a linear vector space. Its dimension is  called  the \emph{degree of mobility} of $g$. Locally,  the {degree of mobility} of $g$ coincides with the dimension of the set (equipped with its  natural topology) of metrics geodesically equivalent to $g$.

The degree of mobility is at least one (since $\const \cdot g$ is always  geodesically equivalent to $g$) and is at most $(n+1)(n+2)/2$, which is the degree of mobility of simply-connected spaces of constant sectional curvature.  
 
Our main result is:

\begin{Th} \label{thm1} 
Let $g$ be a complete Riemannian or pseudo-Riemannian metric on a connected  $M^{n}$ of dimension  $n\ge 3$.  Assume that for every constant $c\ne 0$ the metric  $ c\cdot g$ is not the  Riemannian metric of  constant curvature $+1$. 

If the degree of mobility of the metric is $\ge 3$, then  every complete metric $\bar g$ geodesically equivalent to $g$ is affine equivalent to $g$. 
\end{Th}

 The assumption that the  metrics are complete is important: 
the examples  constructed by Solodovnikov \cite{s1,s3}     show the existence of complete metrics with big degree of mobility (all metrics geodesically equivalent to such metrics are not complete).

\begin{Th} \label{thm2} 
Let $g$ be a  complete Riemannian or pseudo-Riemannian metric on a   closed (=compact, without boundary) connected manifold   $M^{n}$ of dimension  $n\ge 3$. Assume that for every constant $c \ne 0$ the metric  $c \cdot g$ is not the  Riemannian metric of  constant curvature $+1$. 
Then, at least one  the following   possibilities holds:  
\begin{itemize} 
\item the degree of mobility of $g$ is at most two,  or 
 \item every metric $g$ geodesically equivalent to $\bar g$ is affine equivalent to $g$.
\end{itemize} 
\end{Th}

\begin{Rem}  In the Riemannian case, Theorem \ref{thm1} was proved in \cite[Theorem 16]{archive} and in \cite{sbornik}.    The proof uses observations which are wrong in the pseudo-Riemannian situation; we comment on   them in Section \ref{lich}. Our proof for the pseudo-Riemannian case is also not applicable in the Riemannian case, since it uses lightlike geodesics in an essential way.  In Section \ref{proof2},  we give a new, shorter  (modulo  results of our paper)  proof of Theorem \ref{thm1} for the Riemannian metrics as well. 
\end{Rem} 

\begin{Rem} In the Riemannian case, Theorem \ref{thm2} follows from Theorem \ref{thm1}, 
since   every  Riemannian  metric  on a closed manifold is complete. In the 
pseudo-Riemannian case,  Theorem \ref{thm2} is a separate statement.  \end{Rem} 

\begin{Rem}  Moreover, the assumptions  that the metric   is complete and the dimension is $\ge 3$   could be removed from Theorem \ref{thm2}: by  \cite[Corollary 5.2]{mounoud} and \cite[Corollary 1]{quadratic_global}, {\it if 
the degree of mobility of  $g$ on closed $(n\ge 2)-$dimensional manifold 
 is at least three, then for a certain constant $c\ne 0$ the metric $c\cdot g$ is the Riemannian metric of curvature $1$, or every metric geodesically equivalent to $g$ is affine equivalent to $g$. } 
 
 The proofs in \cite{mounoud} and \cite{quadratic_global}  are nontrivial; the proof of  \cite[Corollary 5.2]{mounoud} is  in particularly based on the results of Section  \ref{Const} of the present paper.

\end{Rem} 

\subsection{Motivation I:  projective Lichnerowicz conjecture} \label{lich}

 Recall that a \emph{projective transformation} of the manifold $(M, g)$ is  a diffeomorphism of the manifold that takes (unparametrized) geodesics to geodesics. The infinitesimal generators of the group of projective transformations are complete \emph{projective vector fields}, i.e., complete  vector fields  whose flows take  (unparameterized) geodesics to geodesics. 

Theorem \ref{thm1} allows us  to prove an important partial case of  
the following conjecture,  which folklore attributes  (see \cite{archive} for discussion)  to Lichnerowicz  and Obata (the latter  assumed in addition that the manifold is closed, see,  for example, \cite{hasegawa,nagano,Yamauchi1}):

\emph{{\bf Projective Lichnerowicz Conjecture. } Let a connected  Lie group $G$  act on a complete
 connected  pseudo-Riemannian manifold $(M^n, g)$ of dimension
 $n\ge 2$   by projective
transformations.  Then, it  acts by affine transformations, or  for a certain $c\in \mathbb{R}\setminus \{0\}$ the metric $c\cdot g$
is the Riemannian  metric of  constant positive sectional  curvature $+1$.
 }

 We see that Theorem   \ref{thm1}  implies  
 
 \begin{Cor} \label{lichcon} The 
projective Lichnerowicz  Conjecture  is true under the additional assumption that the dimension $n\ge 3$ and that the degree of mobility of the metric $g$ is $\ge 3$. 
  \end{Cor}

 Indeed,  the pullback of the (complete) metric $g$ under the   projective transformation is a complete  metric geodesically equivalent to $g$. 
 Then, by Theorem \ref{thm1},  it  is affine equivalent to $g$, i.e., 
 the projective transformation is actually  an affine transformation, as it is stated in Corollary  \ref{lichcon}.

 Corollary \ref{lichcon} is  thought  to be  the most complicated part of  the solution of the projective Lichnerowicz  conjecture for pseudo-Riemannian metrics. We do not know yet whether the Lichnerowicz conjecture is true (for pseudo-Riemannian metrics), but we expect  that its solution  (= proof or   counterexample)  will require no new additional 
ideas beyond those from the Riemannian case.

  To support this optimistic expectation, let us recall   that  the  projective Lichnerowicz conjecture was recently proved for Riemannian metrics \cite{CMH,archive}.  The proof contained three  parts:
 
 \begin{itemize} 
\item[(i)]  proof for the metrics with the degree of mobility $ 2$ (\cite[Theorem 15]{archive}, \cite[Chapter 4]{CMH}),  
\item[(ii)] proof under the assumption  $dim(M)\ge 3$ for the metrics with the degree of mobility $\ge 3$ (\cite[Theorem 16]{archive}),  
\item[(iii)]  proof under the assumption  $dim(M)= 2$ for the metrics with the degree of mobility $\ge 3$, \cite[Corollary 5 and  Theorem 7]{CMH}.  
 \end{itemize}

  The most complicated (=lengthy; it is spread over \cite[ \S\S 3.2--3.5, 4.2]{archive}) part was  the proof under the additional assumptions (ii). 
 
 The proof was based on the Levi-Civita description of geodesically equivalent metrics, on the calculation  of curvature tensor for  Levi-Civita metrics with  degree of mobility $\ge 3$ due to Solodovnikov \cite{s1,s3}, and on global ordering of  eigenvalues of $a_{i}^j:= a_{ip} g^{p j} $, where $a_{ij} $ is a solution of \eqref{basic},  due to \cite{benenti,zametki,dedicata}. This  proof  can not   be generalized to  the pseudo-Riemannian metrics.  More precisely, a pseudo-Riemannian analog of the Levi-Civita theorem is much more complicated,  calculations of Solodovnikov essentially use positive-definiteness of the metric, and, as examples show,  the global ordering of  eigenvalues of $a_{i}^j$ is simply wrong for pseudo-Riemannian metrics. 
 
  Thus, Theorem  \ref{thm1}  and  Corollary \ref{lichcon} close the a priori most difficult  part of the solution  of the Lichnerowicz conjecture  for the pseudo-Riemannian metrics.

 Let us now comment on (i), (iii), from the viewpoint of the possible generalization of the Riemannian proof  to  the  pseudo-Riemannian case. 
 We  expect   that this is possible. More precisely, the proof of (i) is based on a trick invented by Fubini \cite{Fubini} and  Solodovnikov  \cite{s1}, see also   \cite{obata,sol,CMH}. The trick uses the assumption that the degree of mobility is  two to double the number of PDEs   (for a vector field $v$ to be  projective  for the metric $g$), and to lower the order of  this equation (the initial equations have  order 2, the equations that we get after applying the trick have order 1). This of courses makes everything  much easier; moreover, in the Riemannian case, one  can explicitly solve this system \cite{Fubini,Obata,s1}. After doing this, one has  to analyze   whether the metrics and the projective field are complete; in the Riemannian case    it was possible to do.

 The trick survives in the pseudo-Riemannian setting.  The obtained system of PDE was solved for the  simplest situations (for small dimensions \cite{bryant,alone}, or under
 the additional assumption that the minimal polynomial of  $a_j^i$ coincides with the characteristic polynomial). We  expect that the other part of the program could be realized for pseudo-Riemannian metrics as well, though of course it will require a lot of work.  
 
 Now let us comment on  the proof under the assumptions (iii): $dim(M)=2$, degree of mobility is $\ge 3$.   The initial proof of \cite{CMH} uses the description of quadratic integrals of geodesic flows of complete Riemannian metrics  due to \cite{Igarashi}. This  description  has no analog  for pseudo-Riemannian metrics.  Fortunately, one actually does not need this description anymore: in   \cite{bryant,alone} a complete  list of  2-dimensional  pseudo-Riemannian metrics admitting  a projective vector field was constructed; the degree of mobility for all these metrics has been  calculated. The metrics that  are interesting for the setting (iii) are the metrics (2a, 2b, 2c) of \cite[Theorem 1]{bryant} and  (3d) of \cite[Theorem 1]{alone}, because all other metrics admitting projective vector fields have constant curvature or degree of mobility equal to $ 2$.  All these metrics are given by relatively simple  formulas using only elementary functions. In order to prove the projective Lichnerowicz conjecture in the setting  (iii), one needs to understand which  metrics  from this list  could be extended to a bigger domain; it does not seem to be too complicated.  For the metrics  (2a, 2b,2c) of \cite[Theorem 1]{bryant} it was already done in \cite{manno}.

 As a concequence of Theorem \ref{thm1},  we obtain  the following simpler version of the Lichnerowicz conjecture. \\
  \begin{Cor} \label{simple} { Let $\textrm{\it Proj}_o$  (respectively, $\textrm{\it Aff}_o$) be the connected  component of the  Lie group of projective transformations  (respectively, affine transformations) of  a complete
 connected  pseudo-Riemannian manifold $(M^n, g)$ of dimension
 $n\ge 3$. Assume that  for no  constant $c\in \mathbb{R}\setminus \{0\}$ the metric $c\cdot g$
is the  Riemannian  metric of  constant positive curvature $+1$.
Then, the codimension   of $\textrm{\it Aff}_o$  in  $\textrm{\it Proj}_o$ is at most one.   
 } \end{Cor}

  Indeed,  it is well known (see, for example \cite{archive}, or more classical sources acknowelged therein) that a vector field is projective if the $(0,2)-$tensor 
  \begin{equation} \label{a-2} a:= L_v g - \tfrac{1}{n+1} \trace({g^{-1} L_vg}) \cdot g \end{equation}  is a solution of \eqref{basic}, where $L_v$ is the Lie derivative with respect to $v$.  Moreover, the projective vector  field is affine, iff the function \eqref{lam}  constructed by $a_{ij}$ given by \eqref{a-2}   is constant. 
  
  Now, let us  take two infinitesimal generator of the Lie group  $\textrm{\it Proj}_o$, i.e., two complete projective vector fields $v$ and $\bar v$. In order to show that the  the codimension   of $\textrm{\it Aff}_o$  in  $\textrm{\it Proj}_o$ is at most one, it is sufficient to show that a linear combination of these vector fields is an affine vector field.  We consider  the solutions $a:= L_v g - \tfrac{1}{n+1} \trace({g^{-1} L_vg}) \cdot g $  and $\bar a:= L_{\bar v} g - \tfrac{1}{n+1} \trace({g^{-1} L_{\bar v}g}) \cdot g $   of \eqref{basic}. 
  
  If 
  $a$, $\bar a$, and   $g$  are linearly independent, the degree of mobility of the metric is $\ge 3$.  Then, Corollary  \ref{lichcon} implies  $\textrm{\it Proj}_o=\textrm{\it Aff}_o$.

  Thus,   $a$, $\bar a, \  g$ are linearly dependent. Since the function $\lambda:= \tfrac{1}{2}g_{p q} g^{p q}$, i.e., the function \eqref{lam} 
  constructed by $a=g$,  is evidently constant, there exists a nontrivial 
   linear combination $\hat a$ of $a,\bar a$ such that the corresponding $\hat \lambda$ given by \eqref{lam} is constant.     Since the mapping $$v\mapsto a:= L_v g - \tfrac{1}{n+1} \trace({g^{-1} L_vg} )\cdot g $$ is linear, the linear combination of $v, \bar v$ with the same coefficients  is an affine vector field. \qed

   \subsection{Motivation II: new methods for investigation of global behavior of geodesically equivalent metrics} 
The theory of  geodesically equivalent metrics  has a long and  fascinating history.  First  non-trivial examples 
 were discovered  by   Lagrange \cite{lagrange}. Geodesically equivalent metrics were studied by Beltrami \cite{Beltrami}, Levi-Civita \cite{Levi-Civita}, Painlev\'e \cite{painleve}  and other classics.   One can find more historical details in the surveys \cite{Aminova2,mikes} and in the introduction to the papers \cite{quantum,integrable,topology,hyperbolic,fomenko,sbornik,archive,dedicata}.

 The success of general relativity made necessary  to study geodesically equivalent pseudo-Riemannian  metrics. The textbooks  \cite{eisenhart1,hallbook,Petrov,Petrov2} on pseudo-Riemannian metrics   have chapters on geodesically equivalent metrics. In the popular paper \cite{Weyl2}, Weyl  stated a  few interesting open   problems on geodesic equivalence of pseudo-Riemannian metrics.    Recent references (on the connection between geodesically equivalent metrics and general  relativity) include  Ehlers et al  \cite{ehlers1,ehlers2},   Hall  and Lonie \cite{hall,hall3,hall4}, Hall \cite{hall1,hall2}.

In many cases,  local statements about Riemannian metrics could be generalised for the pseudo-Riemannian setting, though sometimes this generalisation is  difficult. As a rule, it is very difficult  to generalize global statements about Riemannian metrics to the pseudo-Riemannian setting. Theory of geodesically equivalent metrics is not an exception: most local results could be generalized. For example, the  most classical results: the  Dini/Levi-Civita  description of geodesically equivalent metrics \cite{Dini, Levi-Civita} and Fubini Theorem \cite{Fubini} were generalised in \cite{Aminova,splitting,appendix,pucacco,kiosak}. 

Up to now, no global (if the manifold is  closed or complete) methods for investigation of geodesically equivalent metrics were generalized for the pseudo-Riemannian setting. More precisely,  virtually every global result on geodesically equivalent Riemannian   metrics was obtained by combining the following methods. 

\begin{itemize} 
\item ``Bochner technique".  This is a group of methods combining local differential geometry and  Stokes theorem.  In the theory of geodesically equivalent metrics, the most standard (de-facto, the only) way  to use Bochner technique was to  use tensor calculus   to 
 canonically obtain a nonconstant  function $f$  such $\Delta_g f= \const\cdot  f$, where $\const \ge 0$, which  of cause can not exist on closed Riemannian manifolds. 

An  example could be derived from  our paper: from the equation \eqref{tanno} it follows, that  $(\Delta_g \lambda )_{,k} = 2(n+1)B \lambda_{,k}$. Thus, for a certain  $\const\in \mathbb{R}$ we have 
 $(\Delta_g (\lambda+ \const) ) = 2(n+1)B (\lambda+\const)$. If $B$ is positive,   $g$ is Riemannian, and $M$ is closed,  this implies that the function $\lambda$ is constant, which is equivalent to the statement that the metrics are affine equivalent. 

  The first application of this technique  in the theory of geodesically equivalent metrics is due to Japan geometry school of Yano, Tanno, and Obata,  see for example \cite{hiramatu}. 
Also, mathematical schools of Odessa and Kazan were extremely strong in this group  of methods, see the review papers \cite{Aminova2,mikes}, and the references inside these papers.   

Of cause, since for pseudo-Riemannian metrics 
the equation $\Delta_g f= \const\cdot  f$ could have solutions for $\const \ge 0$, this technique completely fails in the pseudo-Riemannian case. 

\item ``Volume and curvature  estimations". For geodesically equivalent metrics $g$ and $\bar g$,  the repametrisation  of geodesics  is controlled by a function $\phi $ given by  \eqref{phi}.  This function also controls the difference between Ricci curvatures of $g$ and $\bar g$. Playing with this, one can  obtain   obstructions  for  the existence of positively definite   geodesically equivalent metrics with  negatively   definite Ricci-curvature  (assuming the manifold is closed, or complete with finite volume).   Recent  references include \cite{kim,shen}.  

This method  essentially uses the positive definiteness of the metrics.

\item``Global ordering of eigenvalues of $a_j^i$". The existence of a metric $\bar g$ geodesically equivalent to $g$ implies the existence of integrals of special form (we  recall one of the integrals in Lemma  \ref{integr}) for the geodesic flow of the metric $g$ \cite{MT,quantum,integrable}. In the Riemannian case, analyse of the   integrals implies   global ordering of the eigenvalues of the tensor 
$a_j^i:= \left(\tfrac{\det(\bar g)}{\det(g)}\right)^{ \tfrac{1}{n+1}}\bar g^{ip} g_{pj}$, where $\bar g^{ip}$  is the tensor dual  to $\bar g_{ij}$, see  \cite{benenti,zametki,dedicata}. Combining it  with the Levi-Civita description of geodesically equivalent metrics, one could describe topology of closed  manifolds admitting geodesically equivalent Riemannian metrics  \cite{kruglikov,dim2,ERA,starrheit,japan,topology,hyperbolic,dubrovin,spt}. 

Though the integrability survives in the pseudo-Riemannian setting \cite{benenti,involutivity}, the global ordering of the eigenvalues is not valid anymore (there exist counterexamples), so this method also is  not applicable to the pseudo-Riemannian metrics.

\end{itemize}

Our proofs (we explain the scheme of the proofs in the beginning of  Section \ref{proofs})  use essentially new methods. We would like to emphasize here once more  that 
the last step of the proof, which uses the local results to obtain global statements,   is based on the  existence of lightlike  geodesics, and, therefore, is   essentially pseudo-Reimannian.

 A similar idea was used in \cite{einstein}, where it was proved that  complete  Einstein metrics are geodesically rigid: every complete 
 metric geodesically equivalent to  a complete Einstein metric is affine equivalent to it.

  We expect further application of these new methods in the theory of geodesically equivalent metrics.

  \subsection{ Additional motivation: superintegrable metrics.}

  Recall that a metric is called \emph{superintegrable}, if
the number of independent integrals of special form is greater than
the dimension of the manifold. Superintegrable systems are nowadays a hot topic in mathematical physics, probably because  almost all exactly solvable systems are superintegrable.  There are different possibilities for
the special form of integrals; de facto the most standard special
form  of the integrals is  the so-called Benenti integrals, which
are essentially the same as geodesically equivalent  metrics, see \cite{Benenti3,benenti,era1}.  Theorem \ref{thm2}   of our  paper shows that  complete   Benenti-superintegrable  metrics of
nonconstant curvature cannot exist on closed manifolds, which was a folklore
conjecture.

  {\bf Acknowledgements.}
We  thank    Deutsche Forschungsgemeinschaft (Priority Program 1154 --- Global Differential Geometry and research training group   1523 --- Quantum and Gravitational Fields)
  and FSU Jena for partial financial support, and Alexei Bolsinov  and Mike Eastwood for useful discussions.  We  also thank Abdelghani Zeghib for finding a misprint in the main theorem in the  preliminary  version of the paper, and Graham Hall for his  grammatical and stylistic suggestions.

 \section{Proof of Theorems 1, 2}  \label{proofs} 
 
  In Section \ref{standard},  we recall standard facts about  geodesically equivalent metrics and fix the notation.    
  In Section \ref{alg},  we will prove  Lemma \ref{linalg},  which is a purely  linear algebraic statement. Given two solutions of the equation \eqref{int1}, it gives  us  the  equation \eqref{vnb}. The coefficients in the equation are a priori functions. 
  We will work with this equation for a while: 
   In Section \ref{numbers}, we prove (Lemma \ref{de})   that (under the
    assumptions of Theorem \ref{thm1}) one of    the coefficients of 
    \eqref{vnb} is  actually a  constant. Later, we will show (Lemma \ref{bcon}) that the metric $g$ determines the constant uniquely.  
    
The equation \eqref{vnb} will be used in Section \ref{ode}. 
The main result of this section   is Corollary \ref{evolution1}. This corollary gives us (under assumptions of Theorem \ref{thm1}) an ODE that must  be fulfilled along every lightlike   geodesic, and that controls the   reparameterization that produces  $g$-geodesics from $\bar g$-geodesics. 
   The ODE is relatively simple and could be solved explicitly (Section \ref{proof1}). Analyzing the solutions,  we will see that the geodesic is complete with respect to  both metrics iff the function controlling the reparametrization of the geodesics is a constant, which  implies  that the metrics are affine equivalent.   This proves Theorem \ref{thm1} provided   lightlike geodesics exist.   As we mentioned in the introduction, Theorem \ref{thm1} was already proved \cite{japan,archive} for Riemannian metrics.   Nevertheless, for self-containedness,  in Section \ref{proof2}   we give  a new  proof  for Riemannian metrics as well, which is much shorter  than the original proof from  \cite{japan,archive}.

   The proof of Theorem \ref{thm2} will be done in Section \ref{thelast}. The idea is similar: we analyze  a certain ODE along  lightlike geodesics (this ODE will easily follow from the equation \eqref{tanno}, which is an easy corollary of the equation \eqref{vnb}), 
   and show that the assumption that the manifold is closed implies that the solution of the  ODE coming from the metric $\bar g$ is constant, which  implies that  $g $ and $\bar g $
   are   geodesically equivalent.

 \subsection{ Standard formulas we will use } \label{standard} 
 We work in tensor notation with the background metric $g$. That means, we sum with respect to repeating indexes, use $g$ for  raising and lowering  indexes (unless we explicitly say otherwise), and use the Levi-Civita connection of $g$ for   covariant differentiation.

 As it was known already to Levi-Civita \cite{Levi-Civita},  two connections $\Gamma= \Gamma_{jk}^i $ and  $\bar \Gamma= \bar \Gamma_{jk}^i $  have the same unparameterized geodesics, if and only if their difference is a pure trace: there exists a $(0,1)$-tensor $\phi  $ 
such    that \begin{equation} \label{c1} 
 \bar \Gamma_{jk}^i  = \Gamma_{jk}^i + \delta_k^i\phi_{j} + \delta_j^i\phi_{k}.    
   \end{equation}

The reparametrizations of the geodesics for $\Gamma$ and $ \bar   \Gamma $
connected by \eqref{c1} are done according to the following rule:  for a parametrized geodesic  $\gamma(\tau)$ of $\bar \Gamma$,  the  curve $\gamma(\tau(t)) $ is a parametrized geodesic of   $ \Gamma$,  if and only if the parameter transformation $\tau(t)$ satisfies the following ODE: 
\begin{equation} \label{umparametrisation} 
\phi_p \dot \gamma^p =  
\frac{1}{2} \frac{d}{dt} \left(\log\left(\left|\frac{d\tau}{dt}\right|\right)\right). 
\end{equation} 
(We denote by $\dot \gamma$ the velocity vector of  $\gamma$ with respect to the parameter $t$,  and assume summation with respect to repeating index $p$.) 

If $\Gamma$  and   $\bar \Gamma$  related by  \eqref{c1}  are Levi-Civita connections of  metrics $g$ and $\bar g$, then one can find explicitly (following Levi-Civita \cite{Levi-Civita}) a function $\phi$ on the manifold such that its differential $\phi_{,i}$  coincides with the covector $\phi_i$: indeed, contracting \eqref{c1}  with respect to  $i$ and $j$, we obtain 
 $\bar \Gamma_{p i}^p   = \Gamma_{p i}^p  + (n+1) \phi_{i}$. 
  On the other hand, for the Levi-Civita 
 connection  $\Gamma$ of a metric $g$  we have  $  \Gamma_{p k}^p  = \tfrac{1}{2} \frac{\partial \log(|det(g)|)}{\partial x_k} $.  Thus, 
 \begin{equation} \label{c1,5}  \phi_{i}= \frac{1}{2(n+1)}  \frac{\partial }{\partial x_i}  \log\left(\left|\frac{\det(\bar g)}{\det( g)}\right|\right)= \phi_{,i} \end{equation} 
  for the function $\phi:M\to \mathbb{R}$ given by 
  \begin{equation} \label{phi}  \phi:= \frac{1}{2(n+1)} \log\left(\left|\frac{\det(\bar g)}{\det( g)}\right|\right). \end{equation}  In particular, the derivative of $\phi_i$ is  symmetric, i.e., $\phi_{i,j}= \phi_{j,i}$.

The formula  \eqref{c1}   implies  that two metrics $g$ and $\bar g$ are geodesically equivalent if and only if   for a certain $\phi_{i}$ (which is, as we explained above, the differential of $\phi$ given by  \eqref{phi}) we have 
\begin{equation}\label{LC}
    \bar g_{ij, k} -  2 \bar g_{ij} \phi_{k}-  \bar g_{ik}\phi_{j} -   \bar g_{jk}\phi_{i}= 0, 
\end{equation} 
where ``comma" denotes the covariant derivative with respect to the connection $\Gamma$. 
Indeed, the left-hand side of this equation is the covariant derivative with respect to  $\bar \Gamma$, and vanishes if and only if  $\bar \Gamma$ is the Levi-Civita connection for $\bar g$.

The equations \eqref{LC} can be linearized by a clever  substitution:  consider 
   $a_{ij}$ and $\lambda_i$ given by 
\begin{eqnarray} \label{a}
a_{ij} &=   &e^{2\phi} \bar g^{p q} g_{p i} g_{q j},\\  \label{lambda} 
\lambda_{i} & = &  -e^{2\phi}\phi_p \bar g^{p q} g_{q i}, \end{eqnarray}
where $\bar {g}^{p q}$ is  the tensor dual to $\bar g_{p q}$:  \ $\bar {g}^{p i}\bar g_{p j}= \delta_j^i$. 
It is an easy exercise to show that the following linear  equations for   the symmetric $(0,2)$-tensor $a_{ij}$ and $(0,1)$-tensor $\lambda_i$ are    equivalent to \eqref{LC}. 
 \begin{equation} \label{basic} 
 a_{ij,k}= \lambda_i g_{jk} + \lambda_j  g_{ik}. 
 \end{equation}

 \begin{Rem} 
 For dimension 2, the substitution  (\ref{a},\ref{lambda})  was already known to R.  Liouville \cite{liouville} and Dini \cite{Dini}, see \cite[Section 2.4]{bryant} for details and a conceptual explanation.  For arbitrary dimension,  the substitution (\ref{a},\ref{lambda}) and the equation \eqref{basic} are due to Sinjukov \cite{Sinjukov}.    The underlying geometry  is explained in \cite{eastwood1,eastwood}.
 \end{Rem}

Note that it is possible to find a  function $\lambda$ whose  differential is precisely 
 the $(0,1)$-tensor $\lambda_i$: indeed, multiplying \eqref{basic}  by $g^{ij}$ and summing with respect to repeating indexes $i,j$ we obtain $(g^{ij}a_{ij})_{,k} = 2  \lambda_k$. Thus,
$\lambda_i$ is the differential of the function 
\begin{equation}\label{lam} 
\lambda:= \tfrac{1}{2} g^{p q }a_{p q}.   
\end{equation}  
In particular, the covariant derivative of $\lambda_i$ is symmetric:  
 $\lambda_{i,j} = \lambda_{j,i}$.

We see that the equations \eqref{basic} are linear. Thus the  space of the solutions is a linear vector space. 
Its dimension is called the \emph{ degree of mobility} of the metric $g$.

We will also need  integrability conditions for the equation \eqref{basic}  (one obtains them  substituting  the derivatives of $a_{ij}$ given by \eqref{basic} in the formula  $a_{ij,lk}- a_{ij,kl}= a_{i p }R^{p}_{jkl} +  a_{p j}R^{p}_{ikl}$, which is true for every $(0,2)-$tensor  $a_{ij}$)   
 
\begin{equation}  a_{i p }R^{p}_{jkl} +  a_{p j}R^{p}_{ikl} =\lambda_{ l,i} g_{jk}+\lambda_{ l,j} g_{ik}-\lambda_{ k,i} g_{jl}-\lambda_{ k,j} g_{il}. \label{int1} \end{equation}  
 
 The integrability condition in this form was  obtained by Sinjukov \cite{Sinjukov}; 
  in equivalent form, it  was known to   Solodovnikov \cite{s1}. 

As a consequence of these integrability conditions, we obtain that every solution $a_{ij}$ of \eqref{basic} must commute with the Ricci tensor $R_{ij}$:  \begin{equation} \label{first}  
a^{p}_{i}R_{p j} = a^{p}_{j}R_{i p }. \end{equation}

To show this, we ``cycle" the equation \eqref{int1}  with respect to $i,k,l$, i.e.,  we  sum it with itself after renaming the indexes according to $(i\mapsto k \mapsto l \mapsto i)$ and 
with  itself after renaming the indexes according to $(i\mapsto l \mapsto k \mapsto i)$. 
The first term at the left-hand side of the equation will disappear because of the Bianchi equality $R^{p}_{ikl} +  R^{p}_{kli} + R^{p}_{lik}=0 $, the right-hand side vanishes completely,  
 and we obtain 

\begin{equation} 
a_{p i} R^p_{jkl} + a_{p k} R^p_{jli}+ a_{p l} R^p_{jik}=0.
\end{equation} 

Multiplying with $g^{jk}$, using symmetries of the curvature tensor,   and summing over the repeating indexes we obtain $a_{p i} R^p_l - a_{p l} R^{p }_i=0 $, i.e.,  \eqref{first}.

\begin{Rem}\label{finite} 
For further use, let us recall that the equations \eqref{basic} are of finite type (they close after two differentiations \cite{eastwood,mikes,Sinjukov}). Since they are linear, and since  in view of \eqref{lam} they could be viewed as  equations on $a_{ij}$ only, linear independence of the solutions on the whole connected manifold implies linear independence of the restriction of the solutions to every neighborhood.  Thus, the assumption that the degree of mobility  of $g$ (on a connected $M$) is $\ge 3$ 
 implies that the degree of mobility of the restriction of $g$ to every neighborhood is also $\ge 3$. 
\end{Rem} 
 
 We will also need the following statement from  \cite{MT, dedicata}. We denote by $ \textrm{co}(a)^i_j$ the classical comatrix (adjugate matrix) of the $(1,1)$-tensor $a_j^i$ viewed as an $n\times n$-matrix. $ \textrm{co}(a)^i_j$  is also a $(1,1)$-tensor.

\begin{Lemma}[\cite{MT,dedicata}] \label{integr}
If the $(0,2)$-tensor  $a_{ij}$ satisfies  \eqref{basic}, then the function \begin{equation}\label{integral} I:TM\to \mathbb{R}, \  \ \ (\underbrace{x}_{\in M},\underbrace{\xi}_{\in T_xM})\mapsto g_{pq}\textrm\ {\textrm{co}(a)}^p_\gamma \xi^\gamma\xi^q\end{equation} 
    is an integral of the geodesic flow of $g$. 
\end{Lemma}  

Recall that a function is an \emph{integral} of the geodesic flow of $g$, if it is constant along the orbits of the geodesic flow of $g$, i.e., if for every parametrized  geodesic $\gamma(t)$ the function $I\left(\gamma(t),\dot\gamma(t)\right)$ does not depend on $t$. 

\begin{Rem} 
If the tensor $a_{ij}$ comes  from a geodesically equivalent metric $\bar g$ by formula \eqref{a}, the integral \eqref{integral}  is   $$I(x, \xi)=  \left|\tfrac{\det(g)}{\det(\bar g)}\right|^{2/(n+1)}  \bar g(\xi,\xi).$$ In this form,  Lemma \ref{integr} was already known to  
 Painlev\'e \cite{painleve}.   
\end{Rem}

 \subsection{An algebraic lemma }  \label{alg}

 \begin{Lemma} \label{linalg}  Assume 
 symmetric $(0,2)$ tensors $a_{ij}$, $A_{ij}$, $\lambda_{ij}$  and $\Lambda_{ij}$ satisfy  
 \begin{equation}\begin{array}{lcl}a_{i p }R^{p}_{jkl} +  a_{p j}R^{p}_{ikl} &=&\lambda_{ li} g_{jk}+\lambda_{ lj} g_{ik}-\lambda_{ ki} g_{jl}-\lambda_{ kj} g_{il} \\  A_{i p }R^{p}_{jkl} +  A_{p j}R^{p}_{ikl} &=&\Lambda_{ li} g_{jk}+\Lambda_{ lj} g_{ik}-\Lambda_{ ki} g_{jl}-\Lambda_{ kj} g_{il} , \end{array} 
 \label{int2}\end{equation}where $g_{ij}$ is a metric and $R^{i}_{jkl}$ is its curvature tensor. Assume $a_{ij}, A_{ij},$ and $g_{ij}$ are 
 linearly independent at the point $p$. Then, at the point, $\lambda_{ij} $ is a linear combination of $a_{ij}$ and $g_{ij}$.      
 \end{Lemma}

 \begin{Rem} We would like to emphasize here that, though the lemma is formulated in the tensor notation,   it is a purely  algebraic statement (in the proof we will not use differentiation, and, as we see, no differential  condition on $a, A$ is  required).   Moreover, we can replace $R^i_{jkl}$
  by  any (1,3)-tensor having the same algebraic  symmetries (with respect to $g$)  as the curvature tensor, so that for example the fact that 
  the first equation of \eqref{int2} coincides with \eqref{int1} will not be used in the proof (but of cause this will be used in the applications of Lemma \ref{linalg}).  The underlying algebraic structure  of the lemma is  explained in the last section of \cite{kiosak}.  \end{Rem} 
 
{\bf Proof. } First observe that the equations (15) are unaffected by replacing
$$a_{ij} \mapsto  a_{ij} + a\cdot  g_{ij}\ , \  \lambda_{ij} \mapsto  \lambda_{ij} + \lambda\cdot g_{ij}\ ,  \ A_{ij} \mapsto  A_{ij} + A\cdot g_{ij}\ , \  \Lambda_{ij} \mapsto  \Lambda_{ij} + \Lambda\cdot g_{ij}$$
for arbitrary $a, \lambda, A, \Lambda\in \mathbb{R}$.
 Therefore we may suppose, without loss of generality, that
$a_{ij} , \lambda_{ij} , A_{ij} , \Lambda_{ij}$ are trace-free, i.e., \begin{equation} \label{tracefree}a_{ij}g^{ij}=   \lambda_{ij}g^{ij} =  A_{ij}g^{ij} = \Lambda_{ij}g^{ij}=0.\end{equation}

Our assumptions become that $a_{ij}$  and $A_{ij} $ are linearly
independent and our aim is to show that $\lambda_{ij} = \const\cdot  a_{ij}$. 

 We  multiply the  first equation of  \eqref{int2} by $A^{l}_{l'}$ and sum over $l$. After renaming $l'\mapsto l$, 
  we obtain   \begin{equation}  a_{i p }R^{p}_{jkq} A^{q}_l +  
a_{p j}R^{p}_{ikq} A^{q}_l =\lambda_{ p i} A^{p}_l   g_{jk}+
\lambda_{ p j}  A^{p}_l  g_{ik}-\lambda_{ ki}   A_{jl} -\lambda_{ kj} A_{il}. \label{int3} \end{equation}

 We use  symmetries  of the Riemann tensor to  obtain   $a^{p }_i R_{p j k q} A^{q}_l = a_i^p R_{q k j p } A^{q}_l  =a_i^p A_{q l} R^{q}_{k j p}$. After substituting this in  \eqref{int3}, we get 
 
\begin{equation}  a_{i}^p A_{q l} R^{q}_{ k jp}  +  a_{j}^p A_{q l} R^{q}_{ k ip}  = \lambda_{ p i} A^{p}_l   g_{jk}+
\lambda_{ p j}  A^{p}_l  g_{ik}-\lambda_{ ki}   A_{jl} -\lambda_{ kj} A_{il}. \label{int4} \end{equation}

Let us now symmetrize  \eqref{int4} by $l,k$
\begin{equation}\begin{array}{ll}  &
a_{i}^p \left(A_{q l} R^{q}_{ k jp}  +   A_{q k} R^{q}_{ l jp} \right)+  
a_{j}^p\left( A_{q k} R^{q}_{ l ip}  +   A_{q l} R^{q}_{ k ip}\right)  \\
= & \lambda_{ p i} A^{p}_l   g_{jk}+
\lambda_{ p j}  A^{p}_l  g_{ik}-\lambda_{ ki}   A_{jl} -\lambda_{ kj} A_{il}+
\lambda_{ p i} A^{p}_k   g_{jl}+
\lambda_{ p j}  A^{p}_k  g_{il}-\lambda_{ li}   A_{jk} -\lambda_{ lj} A_{ik}.\end{array} \label{int5}\end{equation}

We see that the components  in brackets are the left-hand side of the  second equation 
of \eqref{int2} with other indexes. Substituting  \eqref{int2} 
 in \eqref{int5}, we obtain 
\begin{equation}\begin{array}{ll}  &
a^p_i \Lambda_{p l} g_{jk} + a_i^p \Lambda_{p k}g_{jl}- \Lambda_{jl} a_{ik} -  \Lambda_{jk} a_{il}+  
a^p_j \Lambda_{p l} g_{ik} + a_j^p \Lambda_{p k}g_{il}- \Lambda_{il} a_{jk} -  \Lambda_{ik} a_{jl}  \\
= &\lambda_{ p i} A^{p}_l   g_{jk}+
\lambda_{ p j}  A^{p}_l  g_{ik}-\lambda_{ ki}   A_{jl} -\lambda_{ kj} A_{il}+
\lambda_{ p i} A^{p}_k   g_{jl}+
\lambda_{ p j}  A^{p}_k  g_{il}-\lambda_{ li}   A_{jk} -\lambda_{ lj} A_{ik}.\end{array} \label{int8}\end{equation}

 Collecting the terms by $g$, we  see that \eqref{int8} is can be written as 
 \begin{equation}\begin{array}{ll}  &
\left(a^p_i \Lambda_{p l} -  \lambda_{ p i} A^{p}_l  \right) g_{jk} 
  + \left(a_i^p \Lambda_{p k}-
\lambda_{ p i} A^{p}_k  \right) g_{jl}
+  \left(a^p_j \Lambda_{p l}  -\lambda_{ p j}  A^{p}_l \right) g_{ik}
+ \left(a_j^p \Lambda_{p k}-
\lambda_{ p j}  A^{p}_k \right) g_{il} \\
=&
 \Lambda_{jl} a_{ik} +  \Lambda_{jk} a_{il} + \Lambda_{il} a_{jk} +  \Lambda_{ik} a_{jl}   -\lambda_{ ki}   A_{jl} -\lambda_{ kj} A_{il}-\lambda_{ li}   A_{jk} -\lambda_{ lj} A_{ik}.\end{array} \label{int9}\end{equation}

After denoting \begin{equation} \label{gover} 
\tau_{il}:= a^{p}_i \Lambda_{ p l} - A^{p}_l \lambda_{ p i}\end{equation} the   equation \eqref{int9}  can be written as  
\begin{equation}\begin{array}{ll}  
&\tau_{il}   g_{jk} 
  + \tau_{ik}  g_{jl}
+  \tau_{jl} g_{ik}
+ \tau_{jk} g_{il} \\
=&
 \Lambda_{jl} a_{ik} +  \Lambda_{jk} a_{il} + \Lambda_{il} a_{jk} +  \Lambda_{ik} a_{jl}   -\lambda_{ ki}   A_{jl} -\lambda_{ kj} A_{il}-\lambda_{ li}   A_{jk} -\lambda_{ lj} A_{ik}.\end{array} \label{int10}\end{equation} 
 Multiplying \eqref{int10} by  $g^{jk}$, contracting with respect to $j,k$,  and using \eqref{tracefree},   we obtain 
   \begin{equation}\begin{array}{rl}  
(n+2) \tau_{il}   
+ \left(\tau_{jk} g^{jk}\right)  g_{il} 
=& \Lambda_{p l} a_{i}^p   +  \Lambda_{ip} a^{p}_l    -\lambda_{p i}   A^{p}_l  -\lambda_{ lp } A_{i}^p \\ \stackrel{\eqref{gover}}{=} & \tau_{il} + \tau_{li}.\end{array} \label{int11} \end{equation} 
We see that the right-hand side is symmetric with respect to  $i,l$. Then, so should be the left-hand-side implying 
$\tau_{il}=\tau_{li}$. Then, the equation \eqref{int11} implies 
$
n\tau_{il}+ \left(\tau_{jk} g^{jk}\right)  g_{il}=0
$
implying $\tau_{il}=0$. Then,  the equation \eqref{int10} reads
 \begin{equation}  
0
=
 \Lambda_{jl} a_{ik} +  \Lambda_{jk} a_{il} + \Lambda_{il} a_{jk} +  \Lambda_{ik} a_{jl}   -\lambda_{ ki}   A_{jl} -\lambda_{ kj} A_{il}-\lambda_{ li}   A_{jk} -\lambda_{ lj} A_{ik}. \label{int10bb}\end{equation}

 We alternate \eqref{int10bb} with respect to $j, k$ to obtain  
 \begin{equation}
0= \Lambda_{jl} a_{ik}   +  \Lambda_{ik} a_{jl}   -\lambda_{ ki}   A_{jl}  -\lambda_{ lj} A_{ik} - \Lambda_{kl} a_{ij}   -  \Lambda_{ij} a_{kl}   +\lambda_{ ji}   A_{kl}  +\lambda_{ lk} A_{ij}     .  \label{int13}\end{equation}

Let us now rename $i\leftrightarrow k$ in \eqref{int13} and add the result with  \eqref{int10bb}. We obtain 
 
$$
 \Lambda_{jl} a_{ik}  +  \Lambda_{ik} a_{jl}     -
  \lambda_{ ki}   A_{jl}      -\lambda_{ lj} A_{ik}=0    . $$ In other words, 
   $\Lambda_{\alpha} a_{\beta}  +  \Lambda_{\beta} a_{\alpha}     =
  \lambda_{ \beta}   A_{\alpha}      +\lambda_{ \alpha} A_{\beta}$,  where $\alpha$ and $\beta$  stand for the symmetric indices $jl$ and  $ik$, respectively.

But it is easy
to check that a non-zero simple symmetric tensor 
$X_{\alpha \beta} = P_\alpha Q_\beta + P_\beta Q_\alpha$  determines its
factors $P_\alpha$  and $Q_\beta$  up to scale and order (it is sufficient to check, for example, by taking $P_\alpha$ and $Q_\beta$ to be
basis vectors). Since $a_{ij}$ and $A_{ij}$ 
are supposed to be linearly independent, it follows that
$\lambda_{ij} = \const \cdot  a_{ij}$, as required. \qed

 \subsection{ Local results }
  Within this section, we assume that $(M,g)$ is a connected  Riemannian or pseudo-Riemannian manifold of dimension $n\ge 3$.  
 Recall that the degree of mobility of a metric $g$  is the dimension of the space of the solutions of \eqref{basic}.

\begin{Lemma} \label{degree}
Suppose that the degree of mobility of $g$ is $\ge 3$. \weg{}{Then for every solution $a_{ij}$ of \eqref{basic}, where $\lambda_i$ is the differential of the function $\lambda$ given by \eqref{lam}, there exists an open dense subset $N$ of $M$ each of whose points admits an open neighborhood $U$, a constant $B$, and a function $\mu$ on $U$, such that the hessian of $\lambda$ satisfies on $U$ the equation}
\begin{equation} \label{vnb}
\lambda_{,ij} =  \mu g_{ij} + B a_{ij}.  \end{equation}
\end{Lemma}


{\bf Proof. }
{If $a = \const\cdot g$, then $\lambda$ is constant and the lemma holds with $N=M$, $\mu\equiv B=0$. Otherwise there exists a solution $A$ of \eqref{basic} such that $a,A,g$ are linearly independent.}
We denote by $\Lambda_i$ the $(0,1)$-tensor {from} equation  \eqref{lambda} corresponding to $A$,  i.e., $\Lambda_i= \Lambda_{,i}$  for $\Lambda:= \tfrac{1}{2} A_{p q} g^{p q}$.

  Then the integrability conditions \eqref{int1}
for the solutions $a$ and  $A$ are given by \eqref{int2}  
(with $\lambda_{ij}= \lambda_{,ij}$ and  $\Lambda_{ij}= \Lambda_{,ij}$).

%
%
%

{Let $N$ be the set of all $x\in M$ which admit a neighborhood on which $a,A,g$ are either pointwise linearly independent or pointwise linearly dependent. Being a union of open sets, $N$ is open. $N$ is also dense in $M$: every nonempty open set $U\subset M$ either consists only of points where $a,A,g$ are linearly dependent, then $U\subset N$; or it contains a point where $a,A,g$ are linearly independent and which is therefore contained in $N$.}

{By definition every point in $N$ has an open connected neighborhood $U$ on which one of two possibilities holds:}

\begin{itemize} \item[(a)] {$a,A,g$ are pointwise linearly independent. Then, by Lemma \ref{linalg}, $\lambda_{,ij} =  \mu g_{ij} + B a_{ij}$, where $\mu$ and $B$ are functions; they are unique and smooth because of linear independence. Our goal is to show that $B$ is actually a constant, this will be done in Section \ref{constant}.}

\item[(b)] {$a,A,g$ are pointwise linearly dependent. Then there exist a nonempty open connected subset $U'$ of $U$ and (smooth) functions $\C{1},\C{2}$ on $U'$ such that on $U'$, we have $a+\C{1}A+\C{2}g\equiv0$ or $A+\C{1}a+\C{2}g\equiv0$. (To see that $\C{1},\C{2}$ can be chosen to be smooth, distinguish three cases: the span of $a,A,g$ has on $U$ pointwise dimension $1$; or $A,g$ are linearly independent somewhere; or $a,g$ are linearly independent somewhere.) We will prove in Section   \ref{numbers} that $\C{1},\C{2}$ are actually constants. (Lemma \ref{de} can be applied here because if $a$ or $A$ had the form $\const\cdot g$ on $U'$, then also on $M$, in contradiction to linear independence.) Thus $a,A,g$ are linearly dependent on $U'$ and therefore on $M$. This contradiction rules out case (b).}

\end{itemize}  

\subsubsection{{Linear} dependence of three solutions over functions implies their linear dependence over numbers.} \label{numbers}
We will use  the following statement (essentially due to Weyl \cite{Weyl}); its proof can  be {found} for example in \cite{dedicata}, see also \cite[Lemma 1 in  Section 2.4]{kiosak}.

\begin{Lemma}\label{weyl} Suppose  $a_{ij}$ and $A_{ij}$ are solutions of \eqref{basic}. Assume
$a = f \cdot A$, where $f $ is a function.  Then $f$ is actually a constant.
\end{Lemma}

Our main goal is the following lemma, which {settles} the case (b) of the proof of Lemma \ref{degree}.

\begin{Lemma} \label{de} Suppose for certain functions $\C{1}, \C{2}$
  the  solutions $a, A$  (of \eqref{basic} on a connected manifold $(M^{n\ge 3}, g)$)
  satisfy  
\begin{equation} \label{3}
a_{ij} = \C{1}g_{ij} + \C{2} A_{ij}.  \end{equation}
We assume in addition that $A$ is not $\const \cdot g$.  
Then{} the functions  $\C{1}, \C{2}$  are  constants.
\end{Lemma}
\begin{Rem} Though we will use that the dimension of the manifold is at least three, the statement is true in dimension two as well provided the curvature of $g$ is not constant, see \cite{kruglikov}.

\end{Rem}
{\bf Proof of Lemma \ref{de}.}  We assume that $\C{1}_{,k}$ or  $\C{2}_{,k}$  {is} not zero {everywhere}, and find a contradiction.

Differentiating  \eqref{3} and substituting  \eqref{basic} and its analog for the solution $A$, we obtain
\begin{equation} \label{a1}
\lambda_{i} g_{jk} + \lambda_j g_{ik} = \C{1}_{,k}g_{ij} + \C{2}\Lambda_i g_{jk} + \C{2} \Lambda_jg_{ik} + \C{2}_{,k}A_{ij},\end{equation}  
which is evidently equivalent to
\begin{equation} \label{a2}
\tau_{i} g_{jk} + \tau_j g_{ik} = \C{1}_{,k}g_{ij} + \C{2}_{,k}A_{ij},\end{equation}
where $\tau_i= \lambda_i - \C{2} \Lambda_i$.
We see that for every fixed $k$ the left-hand side   is a symmetric matrix  of the form $\tau_i v_j + \tau_j v_i.  $   If  $\C{1}_{,k}$ is not proportional to $\C{2}_{,k}$ {at some point $x\in M$}, this will imply that  $g_{ij}$ also is of the form
$\tau_i v_j + \tau_j v_i$ {at $x$}, which contradicts the nondegeneracy of $g$.  {Thus there exists a function $f$ with}  \begin{equation}  \C{1}_{,k}= f\cdot\C{2}_{,k}.\label{a7}\end{equation} {At each point $x$ there exists a nonzero  vector} $\xi=(\xi^k)\in T_xM$ such that $\xi^k\C{2}_{,k}=0$. Multiplying  \eqref{a2} with $\xi^k$ and summing with respect to $k$, we see that the right-hand side vanishe\weg{d}{s}, and obtain the equation
$\tau_i v_j + \tau_j v_i= 0$, where $v_i:= \xi^kg_{ik}$.  Since $v_i\ne 0$, we obtain $\tau_i=0$ {at $x$}; hence the equation \eqref{a2}
reads  $
f\cdot \C{2}_{,k}g_{ij} = -  \C{2}_{,k}A_{ij}$ {everywhere on $M$}. Since  {the covector field $\C{2}_{,k}$ is pointwise nonzero on some nonempty connected open subset $U$ of $M$}, this equation implies
$f\cdot g_{ij}=-A_{ij}$ {on $U$}.
By  Lemma \ref{weyl}, $f$ is {constant on $U$. By Remark 4, it is constant globally,} which contradicts the assumptions.\qed    

\subsubsection{ In dimension 3, only metrics of constant curvature can have the degree of mobility $\ge 3$.} 

\begin{Lemma} \label{conformal} 
 {Assume that} the conformal Weyl tensor $C_{ijk}^h$ of the metric $g$ on ({a }connected) $M^{n\ge 3}$ vanishes.  If the curvature of the metric is not constant, the degree of  mobility of $g$ is at most two.
 \end{Lemma} 
 
 Since the conformal Weyl tensor $C_{ijk}^h$ of any  metric on a   3-dimensional  manifold vanishes, a special case of Lemma \ref{conformal} is 
\begin{Cor} \label{dim3}
{The degree} of  mobility of {each metric} $g$  of nonconstant curvature on $M^3$ is at most two.
\end{Cor}

 {\bf Proof of Lemma \ref{conformal}.} It is well-known that the curvature tensor of  spaces with $C_{ijk}^h=0$ has the form 
 \begin{equation} \label{decomp} 
 R_{ijk}^h= P^h_k g_{ij} - P_j ^h g_{ik} + \delta_k^h P_{ij}- \delta_j^h P_{ik},
 \end{equation} 
 where $P_{ij}:= \tfrac{1}{n-2}\left(R_{ij}- \tfrac{R}{2(n-1)} g_{ij}\right)$ (and therefore 
 $P_k^h = P_{p k} g^{p h}$).  
 We denote by $P$ the trace of $P_k^h$; easy calculations give us $P= \tfrac{R}{2(n-1)}$.

 Substituting  the equations  \eqref{decomp} in the integrability conditions \eqref{int1}, we obtain 
 \begin{equation} \label{44}
 \begin{array}{cl} &a_{p i}P^p_l g_{jk}- a_{p i}P^p_k g_{jl} + a_{li} P_{jk} -a_{ki} P_{jl}
 +a_{p j}P^p_l g_{ik}- a_{p j}P^p_k g_{il} + a_{lj} P_{ik} -a_{kj} P_{il}    \\=&\lambda_{ l,i} g_{jk}+\lambda_{ l,j} g_{ik}-\lambda_{ k,i} g_{jl}-\lambda_{ k,j} g_{il}.\end{array} 
 \end{equation} 
 
 Multiplying \eqref{44} with $g^{jk}$ and summing with respect to repeating indexes, and using the symmetry of $P_{ij}$ due to  \eqref{first},     we obtain 
 
  \begin{equation} \label{45}
 a_{p i}P^p_l    =\lambda_{ l,i} - \tfrac{P}{n}a_{li }+ \tfrac{\hat P}{n}g_{li } + \tfrac{2 \lambda}{n} P_{il}, 
 \end{equation} 
 where $\hat P= g^{q \gamma}a_{p q} P_\gamma^p- \lambda_{ \ ,p}^p.$ 
 Substituting \eqref{45} in \eqref{44}, we obtain

\begin{equation} \label{46} \begin{array}{cl}
{0 =}& \tfrac{2\lambda}{n} P_{il} g_{jk}    - \tfrac{2\lambda}{n} P_{ik} g_{jl} +
\tfrac{2\lambda}{n} P_{jl} g_{ik}    - \tfrac{2\lambda}{n} P_{jk} g_{{i}l} \\[1ex]
&{} +a_{li} P_{jk} -a_{ki} P_{jl}+ a_{lj} P_{ik} -a_{kj} P_{il} -  \tfrac{P}{n} a_{il} g_{jk}    + \tfrac{P}{n} a_{ik} g_{jl} -
\tfrac{P}{n} a_{jl} g_{ik}    + \tfrac{P}{n} a_{jk} g_{il}. \end{array}
\end{equation}

Alternating the equation \eqref{46} with respect to $j,k$, renaming $i \longleftrightarrow k$, and adding the result to  \eqref{46}, we obtain

\begin{equation} \label{47bis}
\tfrac{2\lambda}{n} P_{il} g_{jk}      - \tfrac{2\lambda}{n} P_{jk} g_{{i}l}  + a_{li} P_{jk}  -a_{kj} P_{il} -  \tfrac{P}{n} a_{il} g_{jk}    + \tfrac{P}{n} a_{jk} g_{il}=0,  
\end{equation}
which is evidently equivalent to

\begin{equation} \label{47}
\tfrac{2\lambda}{n} P_{il} g_{jk}      - \tfrac{2\lambda}{n} P_{jk} g_{{i}l}  
+ a_{li}\left(  P_{jk}  -  \tfrac{P}{n}  g_{jk}\right)  -a_{kj} \left(P_{il}    - \tfrac{P}{n}  g_{il}\right)=0{.}   
\end{equation}
{Hence (in view of $P_{jk} -\tfrac{P}{n}g_{jk} \ne 0$ because by assumption the curvature of $g$ is not constant) there exists a nonempty open set $U$ such that every solution $a_{ij}$ of \eqref{basic} is on $U$ a smooth linear combination of $g_{ij}$ and $P_{ij}$.} {Thus} every three solutions {$g$, $a$, $\hat a$} of \eqref{basic} are {on $U$ } linearly {dependent} over functions. By Lemma \ref{de}, they are {on $U$, and therefore everywhere, }linearly dependent over numbers{.}  \qed

 \subsubsection{Case (a) of Lemma 3: proof that $B = \const$}  \label{constant}

 We consider  a neighborhood $U\subseteq  M^{n\ge 3}$ such that $a, A, g$ are linearly independent at every point of the neighborhood; by Lemma \ref{de}, almost every point has such neighborhood. 
 
 \begin{Rem}
 Within the whole paper we understand ``almost everywhere" and ``almost every" in the topological sense: a condition is fulfilled everywhere (or in almost every point) if and only if it the set of the points where it is fulfilled is everywhere dense.
 \end{Rem}
 
  In  the beginning of the proof of Lemma \ref{degree}, we explained that at every point of the neighborhood the equation \eqref{vnb} holds for certain smooth  functions  $\mu$ and $B$.  
     Our goal is to show that $B$  is actually  a constant (on $U$). 
 
 Because of Corollary \ref{dim3}, we can assume $n=\textrm{dim}(M)\ge 4$. Indeed, otherwise 
 by Corollary \ref{dim3} the curvature of the metrics is constant, and the metric is Einstein.  Then, by \cite[Corollary 1]{einstein}, the equation \eqref{vnb} holds.

 Within the proof, we will 
  use the following equations, the first one is \eqref{basic}, the second follows from Lemma \ref{degree}.  
 \begin{equation} \label{2} \left\{\begin{array}{c} a_{ij,k}= \lambda_i g_{jk} + \lambda_j  g_{ik}\\
 \lambda_{,ij} = \mu g_{ij} + B a_{ij}.\end{array}\right.\end{equation} Our goal will be to show that $B $ is constant.  We assume that it is not the case and show that 
 for a certain covector field $u_i$  and functions $\alpha, \beta $ on the manifold we have $a_{ij}= \alpha  g_{ij} + \beta u_iu_j$. Later we will show that this gives a contradiction with the assumption that the degree of mobility is three.

 We consider the equation $\lambda_{i,j} = \mu g_{ij} + B a_{ij}$. Taking the covariant derivative  $\nabla_k$, we obtain 
 
 \begin{equation} \label{lambdaijk} 
\lambda_{i,jk} =  \mu_{,k} g_{ij} + B_{,k} a_{ij} + B a_{ij,k}  \stackrel{\eqref{basic}}{= }\mu_{,k} g_{ij} + B_{,k} a_{ij}+ B\lambda_i g_{jk}
+ B\lambda_j g_{ik}.  \end{equation} 
  By definition of the Riemannian  curvature, we have 
 $\lambda_{i,jk}- \lambda_{i,kj} = \lambda_p R^{p}_{ijk} $. Substituting \eqref{lambdaijk}  in this  equation, we obtain

 \begin{equation}\label{30} 
 \lambda_p  R^{p}_{ijk} = \mu_{,k} g_{ij}  + B_{,k} a_{ij} - \mu_{,j} g_{ik} - B_{,j} a_{ik} + B\lambda_j g_{ik}- B\lambda_kg_{ij}.  
\end{equation} 
 
 Now, substituting the second equation of \eqref{2} in \eqref{int1}, we obtain 
 \begin{equation} \label{comp} 
 a_{p i} R^{p}_{jkl} + a_{p j} R^{p}_{ikl}= B\left( a_{li}  g_{jk} + a_{lj} g_{ik} - a_{ki}g_{jl} - a_{kj}g_{il}\right).\end{equation}

 We multiply this equation by $\lambda^l$  and sum over $l$. Using that $a_{p i} R^{p}_{jk q }\lambda^q $ is evidently equal to
  $a_{i}^p  R^{q}_{kj  p}\lambda_q$,   we obtain  
 
 \begin{equation} \label{comp1} 
 a_{i}^p  R^{q}_{kj  p}\lambda_q    + a_{j}^p  R^{q}_{ki  p}\lambda_q= B\left( a_{iq }  \lambda^q  g_{jk} + a_{jq } \lambda^q g_{ik} - a_{ki}\lambda_j - a_{kj}\lambda_i\right).\end{equation} 
 
 Substituting   the expressions for  $R^{q}_{kj  p}\lambda_q$  and $R^{q}_{ki  p}\lambda_q$, we obtain

\begin{equation} \label{i3} \s{1}_ia_{jk} + \s{1}_ja_{ik} + \s{2}_ig_{jk} + \s{2}_j g_{ki} - B_{,j} a^p_i a_{p k} - B_{,i}
a^p_j a_{p k}  =0,  \end{equation} 
where $\s{1}_i:= a_i^p B_{,p}-\mu_{,i} + 2 B \lambda_i$ and $  \s{2}_i:= a_i^p \mu_{,p} - 2B \lambda_p a^{p}_i$.

Now let us work with \eqref{i3}: we alternate the equation with respect to $i,k$ to obtain: 

    \begin{equation} \label{i4} \s{1}_ia_{jk}  + \s{2}_ig_{jk}   - B_{,i}
a^p_j a_{p k}  -\s{1}_ka_{ji}  - \s{2}_kg_{ji}   + B_{,k}
a^p_j a_{p i}    =0.   \end{equation}

We rename $j \leftrightarrow k$ and add the result to \eqref{i3}: we obtain 

 \begin{equation} \label{i5} \s{1}_ia_{jk}  + \s{2}_ig_{jk}   =  B_{,i}
a^p_j a_{p k}  .  \end{equation} 

\begin{Rem}  \label{rem8} If $B=\const$ on $U$, then     $\s{1}_ia_{jk}  + \s{2}_ig_{jk}=0$.  Since  by  Lemma \ref{weyl} $a_{jk}$ is not proportional to $g_{jk}$, we have $\s{1}_i = 0$, which implies that   $\mu_{,i}= 2 B \lambda_i$.  
\end{Rem}

The condition \eqref{i5} implies that under the assumption  $B\ne \const$  the covectors $\s{1}_i$, $\s{2}_i$ and $B_{,i}$ are collinear: Moreover, for    
for certain functions $ \C{1}, \ \C{2}$  \begin{equation} \label{i6} \C{1}  B_{,i}=\s{1}_i , \ \  \  \C{2}  B_{,i}=\s{2}_i,  \   \ \  \  \ \C{1}a_{jk}  + \C{2}g_{jk}   =  
a^p_j a_{p k}  .  \end{equation}   

Taking the $\nabla_k$ derivative of  the last  formula of \eqref{i6}, we obtain 

$$
\lambda_p a_j^p g_{ik} + \lambda_i a_{jk} + \lambda_p a_i^p g_{jk} + \lambda_j a_{ik}= 
\C{1}_{,k} a_{ij} + \C{2}_{,k}g_{ij}+ \C{1}\lambda_ig_{jk} + \C{1} \lambda_jg_{ik}.$$

Alternating the last formula with respect to $i$ and $k$,  we obtain: 
 
 \begin{equation} \label{i7}
\s{3}_i a_{jk} - \s{3}_k a_{ij} + \s{4}_i g_{jk} -\s{4}_kg_{ij} = 0,\end{equation}  
where $\s{3}_i =  \lambda_i+ \C{1}_{,i}$,  $\s{4}_i =  \lambda_p a^{p}{_i}-  \C{1}\lambda_i+ \C{2}_{,i}$.   Let us explain that this equation imply   either $a_{ij}= \alpha  g_{ij} +\beta u_iu_j$ (which was our goal), or $\s{3} = \s{4}=0$.  

We fix a point $x\in U$ and assume that $\s{3}_i\ne 0$ at the point. Then, $\s{4}_i\ne 0$ as well. 
For every vector $ \xi \in T_xM$  we multiply  \eqref{i7} by $\xi^j $ and sum with respect to $j$. Denoting $A(\xi)_k:= a_{jk}\xi^j$ and $G(\xi)_k:= g_{jk}\xi^j$, we obtain 
\begin{equation}\label{I7}
\s{3}_i A(\xi)_{k} - \s{3}_k A(\xi)_i + \s{4}_i G(\xi)_{i} -\s{4}_kG(\xi)_{i} = 0.
\end{equation}
Then, the (at most two-dimensional) subspaces of $T^*_xM$ generated by $\{\s{3}_i, A(\xi)_{i}\}$ and by  $\{\s{4}_i,  G(\xi)_{i}\}$ coincide. Since the metric $g$ is nondegenerate, 
varying $\xi$ we obtain all possible elements  of $T^*_xM$  
as $G(\xi)_{i}$, so   the  subspaces generated by $\{\s{4}_i, G(\xi)_{i}\}$ are   all possible at most two-dimensional  subspaces   containing $\s{4}_i$, and the subspace generated by   $\{\s{4}_i\}$
 is the intersection of all such subspaces. Similarly,  the subspace generated by   $\{\s{3}_i\}$  is the intersection of subspaces generated by  $\{\s{3}_i, A(\xi)_{i}\}$.   Thus, $\s{3}_i = -\alpha \s{4}_i$ for a certain constant  $\alpha$, and the equation \eqref{i7} looks 
 \begin{equation}\label{I8}
 \s{3}_i ( a_{jk} -  \alpha g_{jk})  - \s{3}_k (a_{ij} - \alpha g_{jk})  = 0.
\end{equation}
We take $\eta\in T_xM$  such that $\eta^k\s{3}_k= 0$,  multiply \eqref{I8} by $\eta^k $ and sum over $k$. We obtain  that $A(\eta)= \alpha G(\eta)$ for all such $\eta$. 
  Thus, for a certain const $\beta$ we have $a_{ij} = \alpha g_{ij} + \beta \s{3}_i \s{3}_j$ as we claimed.  
  
  In the case {where $\s{3}$ and $\s{4}$ vanish identically on $U'$}, using \weg{}{\eqref{i6}, \eqref{basic} and }\weg{}{the }definition of $\s{3}$ and  $\s{4}$, we obtain
$\lambda_\alpha a_i^\alpha = \tfrac{(n+2)\C{1} - \weg{}{2}\lambda }{n+\weg{1}{4}} \lambda_i$, i.e., that $\lambda_\alpha$ is an eigenvector of $a_i^j$.  Differentiating this equation and substituting {\eqref{2}, \eqref{i6}, \eqref{basic}, and $\s{3}=0$},  we obtain {}
\[
\weg{}{\Big(\mu+\C{1}B-\tfrac{(n+2)\C{1}-2\lambda}{n+4}B\Big)a_{ij} = \Big(\tfrac{(n+2)\C{1}-2\lambda}{n+4}\mu -\lambda^p\lambda_p -\C{2}B\Big)g_{ij} -2\lambda_i\lambda_j .}
\]
 {Assume that the coefficient of $a_{ij}$ vanishes identically on $U'$. Since $g_{ij}$ has rank $\geq4$ and $\lambda_i\lambda_j$ has rank $\leq1$, the coefficient of $g_{ij}$ vanishes identically on $U'$, and thus the covector field $\lambda_i$ vanishes identically on $U'$. Differentiating $\lambda_i=0$, and using $\lambda_{ij}=\mu g_{ij}+Ba_{ij}$ and Lemma 5, we see that either $a=\const\cdot g$ on $U'$ and therefore everywhere, in contradiction to our linear independence assumption; or $B\equiv0$ on $U'$, in contradiction to the choice of $U'$. This shows that also in the case $\s{3} = \s{4} \equiv 0$ there exist a nonempty open subset $U''$ of $U'$ and functions $\alpha,\beta$ on $U''$ and a covector field $u$ on $U''$ with $a_{ij}= \alpha g_{ij} + \beta u_iu_j$.}

  Let us now explain that  if  $a_{ij} $ is not proportional to $g$ and 
  $a_{ij} = \alpha(x) g_{ij} + \beta(x) u_i u_j$ for every point $x$ of some neighborhood, 
   then $\alpha$  is a smooth function, and  $\beta$  (resp.  $u_i$) can be chosen to be smooth function (resp.  smooth covector field), probably in a smaller neighborhood. 
   Indeed, under these assumptions $\alpha$ is  the eigenvalue of $a_{i}^j$  of (algebraic and geometric) multiplicity precisely $n-1$. Then, it is a smooth function. Then, $\beta u_i u_j$  is a smooth $(0,2)$-tensor field. 
   Since $a_{ij}$ are $g_{ij}$ are not proportionaly,  $\beta u_i u_j$ is not zero and we can chose $\beta=\pm 1$. Then, we have precisely  two choices for the covector $u_i(x)$ 
    at every point $x$ and in a small  neighborhood we can choose $u_i(x)$ smoothly.

Thus, under the assumptions of  this section, 
for every solution $a_{ij}$ of \eqref{basic},  we have  (for certain functions $\alpha_1, \alpha_2$  and a covector field $u_i$) \begin{equation}
a_{ij}= \alpha_1 g_{ij} + \alpha_2 u_iu_j.\label{A1}\end{equation}

For the solution $A_{ij}$  an  analog of the equation \eqref{A1} holds so  (in a possible smaller neighborhood)  we  also have (for certain functions $\beta_1, \beta_2$  and a covector field $v_i$)
\begin{equation}
A_{ij}= \beta_1 g_{ij} + \beta_2 v_iv_j.\label{A2}
\end{equation} 

Without loss of generality, we can assume that $a_{ij}+ A_{ij}$ (which is certainly a solution of \eqref{basic}) is also not proportional to $g_{ij}$, otherwise we replace $A_{ij}$ by $\tfrac{1}{2}A_{ij}$. Then, 

\begin{equation}
a_{ij} + A_{ij}= \gamma_1 g_{ij} + \gamma_2 w_iw_j.\label{A3}
\end{equation}

Subtracting \eqref{A3} from the sum of \eqref{A1} and \eqref{A2}, we obtain 
\begin{equation} \label{A4}
(\gamma_1- \alpha_1 - \beta_1 )g_{ij} = \alpha_2 u_iu_j+ \beta_2  v_i v_j - \gamma_2 w_iw_j. 
\end{equation} 
Since the tensor $g_{ij}$ is nondegenerate, its rank coincides with 
 the dimension of $M$ that  is at least  $4$. The rank of the  tensor 
$\alpha_2 u_iu_j+ \beta_2  v_i v_j - \gamma_2 w_iw_j $ is   at most three. Thus  
the coefficient $(\gamma_1- \alpha_1 - \beta_1 )$ must vanish, which  implies    that 
\begin{equation} \label{A6}
\alpha_2 u_iu_j+ \beta_2  v_i v_j = \gamma_2 w_iw_j.
\end{equation} 

 We see that the rank of $\alpha_2 u_iu_j+ \beta_2  v_i v_j$ is at most one, which   implies that  $u_i$ is proportional to $v_i$ (the coefficient of the proportionality is a function). 
 Thus  \eqref{A6} implies that $w_i$ is proportional to $u_i$ as well. Thus $a_{ij},$ $ A_{ij} ,$ and $g_{ij}$ are linearly dependent over functions, which  implies  by Lemma  \ref{de} that they are linearly dependent over numbers.  This is a  contradiction to the assumptions, which proves the remaining part of Lemma \ref{degree}.

\subsubsection{The constant $B$ is universal }  \label{2.5.1} 

Let  $(M^{n\ge 3},g)$ be a  connected  pseudo-Riemannian manifold. Assume   the degree of mobility of  $g$ is $\ge 3$, let $(a_{ij}, \lambda_i) $ be a   solution of the equations  \eqref{basic} such that $a_{ij} \ne \const \cdot g_{ij} $ for every $\const \in \mathbb{R}$. 
    Then, in a neighborhood of almost every point there exist a constant $B$ and a function $\mu$ such that the equations \eqref{2} hold.  Note that the constant $B$ determines the function $\mu$: 
    indeed, multiplying \eqref{vnb} by $g^{ij} $ and summing with respect to $i,j$ we obtain 
    $\lambda^i_{\ \ , i} = n\mu - 2B\lambda$. 
    
   Our goal is to prove  the statement announced in the title of the section: we would like to show that the constant $B$ is the same in all such neighborhoods (which in particular implies that the equations \eqref{2} hold at all  points with one universal constant $B$ 
    and one universal function $\mu$). 
   We will need the following 
   
   \begin{Cor} \label{ta}
Let $a_{ij}, \lambda_{i}$ satisfy the equations (\ref{2}) in a neighborhood $U\subseteq  (M, g)$ with a certain constant $B$ and a smooth function $\mu$. Then the function $\lambda$ given by \eqref{lam} satisfies  the equation  \begin{equation}\label{tanno}
\lambda_{,ijk} - B\left(2 \lambda_{,k} g_{ij} + \lambda_{,j} g_{ik} + \lambda_{,i}g_{jk}\right)=0,
\end{equation}
\end{Cor} 

\begin{Rem}

This equation   is a famous one; it  naturally appeared in different parts  of differential geometry.  
   Obata and Tanno  used this equation trying to understand the connection between the eigenvalues of the laplacian $\Delta_g$  and the geometry and topology of the manifold. 
   They observed \cite{Obata,Tanno} that the eigenfunctions corresponding to the second eigenvalue of the Laplacian of the metrics of constant positive curvature $-B$ on the sphere satisfy the equation \eqref{tanno}. 
   
    Tanno  \cite{Tanno} and  Hiramatu \cite{hiramatu} 
    related  the equations to {projective vector fields}. Tanno has shown that for every solution $\lambda $ of this equation     the vector field $\lambda_{,}^{\ i}$ is a projective vector field (assuming $B\ne 0$), Hiramatu proved the reciprocal statement under certain additional assumptions. 
   
    As it was shown by Gallot \cite{G},  see also \cite{Leist,cones,mounoud},   decomposability 
     of the holonomy group of the cone over a manifold  implies the existence of  a nonconstant solution of the equation \eqref{tanno} on the manifold. 
\end{Rem} 

{\bf Proof of Corollary \ref{ta}.} Covariantly differentiating \eqref{vnb} and replacing the covariant derivative of $a_{ij}$ by \eqref{basic} we obtain \eqref{tanno} { from Remark \ref{rem8} if $a\neq\const\cdot g$. If $a=\const\cdot g$, we have $\lambda_{,i}=0$, thus \eqref{tanno} holds as well.} \qed  
   
   \begin{Cor} \label{100} 
     Let the degree of mobility of a metric $g$ on a connected $(n>3)$-dimensional $M$ be $\ge 3$. Assume $(a_{ij}, \lambda_i)$ is a solution of \eqref{basic}. Then, if  $\lambda_i\ne 0$ at a point, then the set of the points such that $\lambda_i\ne 0$ is everywhere dense. 
   \end{Cor} 
    \begin{Rem} 
 The assumption that the degree of mobility of $g $ is  $\ge 3$   is important: Levi-Civita's description of geodesically equivalent metrics \cite{Levi-Civita}   immediately gives  
 counterexamples.  
 \end{Rem} 
 {\bf  Proof of Corollary \ref{100}.} Combinig Lemma \ref{degree}, Remark \ref{rem8},  and  Corollary \ref{ta}, we 
  obtain that in a neighborhood of almost every point $\lambda$ given by \eqref{lambda} satisfies \eqref{tanno}.  By \cite[Proposition 2.1]{Tanno},   the vector field $\lambda^i$ is a projective vector field  (almost everywhere, and, therefore, everywhere) on $(M, g)$. As it was shown for example in  \cite[Theorem 21.1(ii)]{hallbook}, if it is not zero at a point, then it is not zero at almost every point. \qed 
 
   \begin{Cor}  \label{temp3} Let $a_{ij}, \lambda_{i}$ satisfy the equations (\ref{2}) in a neighborhood $U$ with a certain constant $B$ and a smooth function $\mu$.  Let 
   $\lambda$ be the function constructed by \eqref{lam}. Then for every  geodesic $\gamma(t)$   the following equation holds (at every $t\in \gamma^{-1}(U)$):   
  \begin{equation}\tfrac{d^3}{dt^3} \lambda(\gamma(t))= 4 Bg(\dot\gamma(t), \dot\gamma(t)) \cdot  \tfrac{d}{dt}  \lambda(\gamma(t)) \label{temp1}
, \end{equation}   where $\dot\gamma$ denotes  the velocity vector of the geodesic $\gamma$, and $ g(\dot\gamma(t), \dot\gamma(t)):= g_{ij} \dot\gamma^i \dot\gamma^j$. 
\end{Cor}
   
   {\bf Proof.}  Multiplying \eqref{tanno} by  $\dot\gamma^i \dot\gamma^j \dot\gamma^k$ and summing with respect to $i,j, k$ we obtain \eqref{temp1}.  \qed

 \begin{Lemma} \label{fuck}  Let $(M^{n\ge 3},g)$ be a connected manifold and $(a_{ij}, \lambda_i)$ 
 be a solution of \eqref{basic}. Assume almost every point has a neighborhood such that in this neighborhood there exists a constant $B$ and a smooth function $\mu$ such that the equation \eqref{vnb} is  fulfilled. Then the constant $B$ is the same in  all such neighborhoods. 
 \end{Lemma}

{\bf Proof.}  It is sufficient to prove this statement locally, in a  sufficiently small neighborhood of arbitrary point. We take a small neighborhood  $U$, two points $p_0, p_1\in U$, and two neighborhoods  $U(p_0)\subset U$, $U(p_1)\subset  U$ of these points. We assume that our neighborhoods are  small enough and that  we can connect every point of $U(p_0)$ with every point of $U(p_1)$ by a unique geodesic lying in $U$. 
We assume that   the equation \eqref{vnb} holds  in  $U(p_i)$ with the constant $B:=B_i$; our goal is to show that  $B_0=B_1$.

Suppose  it is not the case. 
 We consider all geodesics  $\gamma_{p,p_0}$ lying in $U$  connecting  all  points $p\in U(p_1)$ with  $p_0$, see the picture. We will think that $\gamma(0)=p_0$  and $\gamma(1)\in U(p_1)$.   
 
 For every such  geodesic $\gamma_{p,p_0}(t)$  there exists a point $q_{p,p_0}:= \gamma_{p,p_0}(t_{p,p_0})$  on this geodesic such that  for all $t\in[0, t_{p,p_0}) $  the following conditions are fulfilled: 
 
 \begin{enumerate} \item 
 the equations \eqref{2} are fulfilled  with $B= B_0$ in a small neighborhood of $ \gamma(t)$, and \item 
 for no neighborhood of $\gamma_{p,p_0}(t_{p,p_0})$ the equations \eqref{2} are fulfilled with $B=B_0$. 
 \end{enumerate} 
 
 Then, at every such point $\gamma_{p,p_0}(t_{p,p_0}) $ we have that $a_{ij}= \tfrac{2}{n} \lambda g_{ij}$. Indeed, the trace-free version of  \eqref{vnb} is 
 \begin{equation} \label{tracefree1}
 \lambda_{,ij}-\tfrac{1}{n}\lambda_{,k}^{\ \ k} = B(a_{ij} - \tfrac{2}{n} \lambda g_{ij})
 \end{equation}
 implying that $B$  is  the coefficient of proportionality of two smooth tensors. If $a_{ij}\ne \tfrac{2}{n} \lambda g_{ij}$ at $\gamma_{p,p_0}(t_{p,p_0})$, we have 
 $ a_{ij} - \tfrac{2}{n} \lambda g_{ij}\ne 0$,  and $B$ can be prolonged to  a smooth function in a small neighborhood of $\gamma_{p,p_0}(t_{p,p_0})$. Since it is locally-constant, it is  (the same)  constant at all points of the  neighborhood of   $\gamma_{p,p_0}(t_{p,p_0}) $ contradicting the  conditions 1,  2. 
 
 Moreover, at every such  point $\gamma_{p,p_0}(t_{p,p_0})$ we have $\lambda_i=0$. Indeed, otherwise we  multiply  \eqref{tanno} by $g^{ij}$ and sum with respect to $i,j$.  We obtain $\lambda^i_{\ \ ,ik}= 2(n+1)B\lambda_k$. We again have that  $B$  is  the coefficient of proportionality of two smooth tensors. Arguing as  above we obtain that $\lambda_i=0$ at every point $\gamma_{p,p_0}(t_{p,p_0})$.

 Since at every point $\gamma_{p,p_0}(t_{p,p_0})$ we have $\lambda_i=0$, we have that 
 $\tfrac{d}{dt} \lambda(\gamma_{p,p_0}(t))_{|t=t_{p,p_0}}= 0$. Then, the set of all such $ \gamma_{p,p_0}(t_{p,p_0})$ contains a smooth (connected)  hypersurface (because the set of zeros of the derivatives of the solutions of the equation \eqref{temp1} depends smoothly on the initial data and on $g(\dot\gamma,\dot\gamma)$).  We denote this hypersurface by $H$. 
 
Since $\lambda_i=0$ at every point of $H$, the function $\lambda $ is constant (we denote it by $\tilde \lambda\in \mathbb{R}$) on $H$.   
\begin{figure} 
 \centerline{{{\psfig{figure=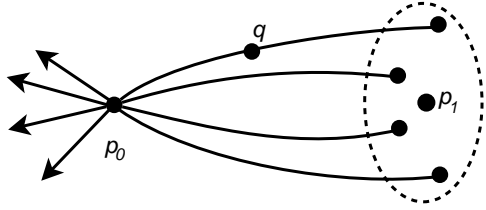,scale=0.4}}}}  
\caption{ The geodesics $\gamma_{p,p_0}$, their velocity vectors at $p_0$, and   the point $q_{p,p_0}=\gamma_{p,p_0}(t_{p, p_0})$ on one of these geodesics }
\end{figure}

 Now let us return to the geodesics  
 $\gamma_{p,p_0}$ connecting    points $p\in U(p_1)$ with $p_0$. 
 We consider  the integral $I$ given by \eqref{integral}. Direct calculations show 
 that at every  point $q$  where  $a_{ij}= c\cdot g_{ij}$ 
  the integral is given by \begin{equation} \label{fk} I(\xi)=c^{n-1}g(\xi, \xi) \end{equation}  (for every tangent vector $\xi\in T_qM$). 
 As we explained above, every such geodesic passing through a point of $H$  has a point 
 such that $a_{ij}= c \cdot g_{ij}$, where $c=\tfrac{2}{n}\tilde \lambda$ is a  constant. Since the integral is constant on the orbits, we have   that   $I\left(\dot\gamma_{p,p_0}(0)\right)= c^{n-1}\cdot g\left(\dot\gamma_{p,p_0}(0), \dot\gamma_{p,p_0}(0)\right)$. 
 Then, the measure of the subset $$\{ \xi\in T_{p_0}M\mid I(\xi)=c^{n-1}\cdot g(\xi, \xi)\}\subseteq T_{p_0}M$$ 
 is not zero.  Since this  set  is given by an  algebraic 
 equation,   it must coincide with the whole $T_{p_0}M$. Then, $a_{ij}= c \cdot g_{ij}$ at 
 the point $p_0$. Since we can replace $p_0$ by every point of its neighborhood $U(p_0)$, we obtain 
 that $a_{ij}= c^{n-1}  \cdot  g_{ij}$ at every point of $U(p_0).$ By  Remark \ref{finite}, \  $a=c^{n-1} \cdot g$ on the whole manifold. 
  \qed

\subsubsection{ The metric $g$ uniquely determines $B$.  } \label{Const} 

By Lemma \ref{degree}, under the assumption that the degree of mobility is $\ge 3$, for every solution $a$ of \eqref{basic} there exists a constant $B$ such that
 the equation \eqref{vnb} holds {on a suitable open set}. In this chapter we show that the constant $B$ is the same for all (nontrivial)  solutions $a_{ij}$, i.e., the metric determines it uniquely. 

\begin{Lemma}  \label{bcon} 
Suppose two nonconstant functions $f, F:M^n\to \mathbb{R}   $ on a connected manifold $(M^n, g)$  of dimension $n>1$ satisfy 
\begin{equation}\label{tanno1}
 \begin{array}{l}
f_{,ijk} - b\left(2 f_{,k} g_{ij} + f_{,j} g_{ik} + f_{,i}g_{jk}\right)=0, \\ 
F_{,ijk} - B\left(2 F_{,k} g_{ij} + F_{,j} g_{ik} + F_{,i}g_{jk}\right)=0,\end{array}
\end{equation} where $b$ and $B$ are constants.  Assume that there exists a point where the
derivative of $f$ is nonzero  and a point where  the derivative of $F$ is  nonzero. Then, $b=B$. 
\end{Lemma} 

{\bf Proof.} By definition of the curvature,  for every function $f$, 
we have   $f_{,ijk}- f_{,ikj}= f_p R^p_{ijk}$; replacing $f_{,ijk}$ by the right-hand side of the first equation of \eqref{tanno1} we obtain. 
\begin{equation}\label{c10}  
f_{,p} R^p _{ijk} = b\left(f_{,k} g_{ij} -f_{,j} g_{ik}\right).
\end{equation} 
The same is true for the second equation of \eqref{tanno1}: 
\begin{equation}\label{c2}  
F_{,p} R^p _{ijk} = B\left(F_{,k} g_{ij} -F_{,j} g_{ik}\right).
\end{equation} 

Multiplying \eqref{c10} by  $F_{,}^{\ k}$ , summing with respect to repeating indexes and using \eqref{c2} we obtain 

\begin{equation}\label{c3}  
B \left(                 F_{,p} f_{,}^{\ p}  g_{ij} - F_{,j} f_{,i} \right)=  
b \left( F_{,p} f_{,}^{\ p}  g_{ij} -F_{,i} f_{,j} \right).
\end{equation} 

Multiplying by $g^{ij} $ and summing with respect to repeating indexes, we obtain 
$
B(n-1)                F_{,p} f_{,}^{\ p} =  
b (n-1)F_{,p} f_{,}^{\ p}   .
$
If $F_{,p} f_{,}^{\ p}\ne 0$ we are done: 
 $B=b$. Assume $  F_{,p} f_{,}^{\ p}=0$. Then,  \eqref{c3} reads  
 $
B                  F_{,j} f_{,i}=  
b F_{,i} f_{,j} .
$ Since by Corollary \ref{100} there exists a point where $F_{,j}$ and $f_{,i}$ are both nonzero, we
obtain again $B = b$.
Then, $f_{,i}$  is proportional to $F_{,j}$. Hence,  $B=b$. \qed  

\subsubsection{ An ODE  along geodesics }    \label{ode}

 \begin{Lemma} \label{const}
 Let $g$ be a metric on a connected $M^{n\ge 3}$ of degree of mobility $\ge 3$. For a metric  $\bar g $  geodesically equivalent to $g$,   let us  consider  $a_{ij}$,   $\lambda_i$, and $\phi$  given by (\ref{a},\ref{lambda},\ref{phi}). Then, the exist  constants $B$, $\bar B$ such that
   the following formula holds:
 \begin{equation}\label{f1} 
 \phi_{i,j} - \phi_{i} \phi_{j} = -B g_{ij}+\bar B \bar g_{ij}  .
 \end{equation} 
 \end{Lemma}

 {\bf Proof. }   
 We covariantly  differentiate \eqref{lambda} (the index of differentiation is ``j"); then we substitute the expression \eqref{LC} for $\bar g_{ij,k}$   to obtain 
 \begin{equation} \label{f2} \begin{array}{ccl}
 \lambda_{i,j} &=& -2 e^{2\phi}\phi_{j} \phi_p \bar g^{p q} g_{q i}-e^{2\phi}\phi_{p,j} \bar g^{p q} g_{q i}+e^{2\phi}\phi_p  \bar g^{p s} \bar g_{s l,j} \bar g^{l q} g_{q i} \\ &=&   -e^{2\phi}\phi_{p,j} \bar g^{p q} g_{q i}+e^{2\phi}\phi_p \phi_s \bar g^{p s}   g_{ i j }+  e^{2\phi}\phi_{j} \phi_l \bar g^{lq}g_{q i}      \end{array}    ,  
 \end{equation} 
 where $\bar g ^{p q}$ is the tensor dual to $\bar g_{p q}$, i.e.,  $\bar g ^{p i}\bar g_{p j}= \delta_j^i$ . 
 We now  substitute  $\lambda_{i,j}$ from \eqref{vnb},  use that $a_{ij}$ is given by \eqref{a}, and divide by $e^{2\phi}$ for cosmetic reasons   to  obtain 
 \begin{equation} \label{f3} 
 e^{-2\phi} \mu g_{ij} + B \bar g^{p q} g_{p j}g_{q i} = -\phi_{p, j} \bar g^{p q} g_{q i}+\phi_p \phi_{s} \bar g^{p s} \bar  g_{ i j }+ \phi_{j} \phi_{l} \bar g^{lq}g_{q i}.  
 \end{equation}  
 Multiplying with $g^{i\xi} \bar g_{\xi k}$,  we obtain 
 \begin{equation} \label{f4} 
 \phi_{k,j}-\phi_{k}\phi_{j} =\underbrace{(\phi_p \phi_q \bar g^{p q} - e^{-2 \phi } \mu )}_{\bar b} \bar g_{kj} - B g_{kj}.  
  \end{equation} 
   The same holds with the roles of $g$ and $\bar g$ exchanged  (the function \eqref{phi} constructed by  the interchanged pair $\bar g, g$ is evidently equal to $-\phi$).   We  obtain 
  \begin{equation} \label{f4b}
   -\phi_{k;j}-\phi_{k}\phi_{j} =\underbrace{(\phi_p \phi_q g^{p q} - e^{2 \phi } \bar \mu )}_{b}  g_{kj} - \bar B \bar g_{kj},
  \end{equation}  
  where $\phi_{i;j}$ denotes the covariant derivative  of $\phi_i$ with respect to the Levi-Civita connection of the metric $\bar g$. 
Since the Levi-Civita connections of $g$ and of $\bar g$ are related by the formula  
(\ref{c1}), we have 
$$
  -\phi_{k;j}-\phi_{k}\phi_{j} = \underbrace{-\phi_{k,j} + 2 \phi_k\phi_j}_{ -\phi_{k;j}} - \phi_{k}\phi_{j}= -(\phi_{k,j}-\phi_{k}\phi_{j}).$$
  We see that the left hand side of \eqref{f4} is equal to minus the left hand side of \eqref{f4b}.
  Thus, $b\cdot g_{ij} - \bar B \cdot \bar g_{ij} = B\cdot g_{ij} - \bar b \cdot \bar g_{ij}$ holds on $U$. Since the metrics $g$ and $\bar g$ are not proportional on $U$ by assumption, $\bar b= \bar B$, and the formula \eqref{f4} coincides with \eqref{f1}.  \qed

\begin{Cor} \label{evolution}
Let $g$, $\bar g$  be geodesically equivalent  metrics on a  connected  $M^{n\ge 3} $ such that 
the degree of mobility of $g$ is $\ge 3$. We consider  a (parametrized)  geodesic $\gamma(t)$ of the metric $g$,   and denote by $\dot \phi$, $\ddot \phi$ and $\dddot \phi$ the first, second and third derivatives of the function $\phi$  given by \eqref{phi} along the geodesic. Then, there exists  a constant $B$ such that   for every geodesic  $\gamma$  
 the following ordinary differential  equation holds: 
   \begin{equation}  
 \begin{array}{lcl}\dddot\phi&=&  4  B g(\dot\gamma, \dot\gamma) \dot\phi  + 6 \dot\phi\ddot\phi- 4(\dot \phi)^3\end{array}, 
  \label{phi0}
 \end{equation} where  $g(\dot\gamma, \dot\gamma):=  g_{ij} \dot\gamma^i\dot\gamma^j$.  
   \end{Cor} 
   Since lightlike geodesics have $g(\dot\gamma, \dot\gamma)=0$ at every point,  a partial case of Corollary \ref{evolution} is 
   \begin{Cor} \label{evolution1}
Let $g$, $\bar g$  be geodesically equivalent  metrics on a  connected  $M^{n\ge 3} $ such that 
the degree of mobility of $g$ is $\ge 3$. Consider a (parametrized) lightlike    geodesic $\gamma(t)$ of the metric $g$, 
 and denote by $\dot \phi$, $\ddot \phi$ and $\dddot \phi$ the first, second and third derivatives of the function $\phi$  given by \eqref{phi} along the geodesic. Then, along the geodesic,   
 the following ordinary differential  equation holds: 
   \begin{equation}  
 \begin{array}{lcl}\dddot\phi&=&   6 \dot\phi\ddot\phi- 4(\dot \phi)^3\end{array}. 
  \label{phi0b}
 \end{equation}
   \end{Cor} 
   
  {\bf Proof of Corollary \ref{evolution}. }  If $\phi\equiv 0$ in a neighborhood $U$,  the equation is automatically fulfilled.  
 Then, it is sufficient to prove Corollary \ref{evolution} assuming $\phi_{i}$ is not constant. 
 
  The formula \eqref{f1}  is evidently equivalent to 
 \begin{equation}
 \phi_{i,j}=\bar B \bar g_{ij}- B g_{ij}+   \phi_{i}\phi_{j}. 
 \label{phi2} 
 \end{equation}
Taking the  covariant derivative of  \eqref{phi2},  we obtain 
  \begin{equation}
 \phi_{i,jk}=\bar B \bar g_{ij,k} +    \phi_{i,k}\phi_{j}+  \phi_{j,k}\phi_{i}. 
  \label{phi3} 
 \end{equation}
 Substituting the expression for $\bar g_{ij,k} $ from \eqref{LC}, and  substituting $\bar B \bar g_{ij}$ given by  
\eqref{f1}, we obtain  
 \begin{equation}
 \begin{array}{lcl}\phi_{i,jk}&=&\bar B  ( 2 \bar g_{ij} \phi_{k}+   \bar g_{ik}\phi_{j} +    \bar g_{jk}\phi_{i})+     \phi_{i,k}\phi_{j}+  \phi_{j,k}\phi_{i} \\&=& B(  2  g_{ij} \phi_{k}+   g_{ik}\phi_{j} +   g_{jk}\phi_{i} ) + 2 (\phi_{k}\phi_{i,j} + \phi_{i}\phi_{j,k}+\phi_{j}\phi_{k,i})- 4\phi_{i}\phi_{j}\phi_{k} \end{array}
  \label{phi4}
 \end{equation}
 Contracting with $\dot\gamma^i \dot\gamma^j \dot\gamma^k$ and using that $\phi_i$ is the differential of the function \eqref{phi}  we obtain the desired ODE \eqref{phi0}. \qed

 \subsection{ Proof of Theorem 1 for pseudo-Riemannian metrics}
  
  Let $g$  be a metric on a connected $M^{n\ge 3}$. Assume that  for no constant $c\ne 0$ the metric $c\cdot g$ is Riemannian,  which in particular implies the   existence of  lightlike geodesics. 
  
    Let $\bar g$ be  geodesically equivalent to $g$. Assume  both metrics are complete. 
 Our goal is to    show that $\phi$ given by \eqref{phi} is constant, because  in view of \eqref{c1}  this  implies that the metrics are affine equivalent.
 
 Consider a parameterized lightlike geodesic $\gamma(t)$ of $g$.  
  Since the metrics are geodesically equivalent, for a certain function $\tau:\mathbb{R}\to \mathbb{ R}$  the 
   curve $\gamma(\tau)$ is a geodesic of $\bar g$.     
 Since the metrics are complete, 
 the reparameterization  $\tau(t)$ is a diffeomorphism $\tau:\mathbb{R}\to \mathbb{R}$. Without loss of generality we can think that $\dot{\tau}:= \tfrac{d}{dt}\tau$ is positive, otherwise we replace 
 $t$ by  $-t$. Then, the equation    \eqref{umparametrisation} along the geodesic reads 
  \begin{equation} 
   \phi(t) = \tfrac{1}{2} \log(\dot\tau(t)) + \const_0 \label{un2} . \end{equation}

 Now let us consider the equation \eqref{phi0b}.  Substituting    
 \begin{equation} 
   \phi(t)= -\tfrac{1}{2}\log(p(t)) + \const_0  \label{un5}  \end{equation} in it (since $\dot\tau>0$, the substitution is global),  we obtain
 \begin{equation} \dddot p = 0 \label{un3} . \end{equation}

  The  
 solution of \eqref{un3} is $p(t) =C_2 t^2 + C_1 t +C_0$.  Combining \eqref{un5} with \eqref{un2}, we see that  $\dot\tau = \tfrac{1}{C_2 t^2 + C_1 t +C_0}$.  Then
   \begin{equation} 
   \tau(t) = \int_{t_0}^t \tfrac{d\xi }{C_2 \xi^2 + C_1 \xi +C_0}\  \ +\const.
  \end{equation} 
  We see that if the polynomial $ C_2 t^2 + C_1 t + C_0$  has real roots (which is always the case if $C_2=0$, $C_1\ne 0$), then  the   integral 
   explodes in finite time. If the polynomial has no real roots, but $C_2\ne 0$,  the function $\tau$ is bounded. Thus, the only possibility for  $\tau $ to  be a diffeomorphism is $C_2=C_1=0$ implying  $\tau(t)  = \tfrac{1}{C_0} t  + \const_1$ implying $\dot\tau =\tfrac{1}{ C_0}$ implying $\phi$ is  constant along the  geodesic.

    Since every two points of a connected  pseudo-Riemannian manifold such that for no constant $c$ the metric $c\cdot g$ is Riemannian  
     can be connected  by a sequence of lightlike geodesics,  $\phi$ is a constant, so that $\phi_i\equiv 0$,       and the metrics are affine equivalent by \eqref{c1}. \qed \label{proof1}
   
\subsection{ Proof of Theorem 1 for Riemannian metrics} 
 \label{proof2}  
   As we already mentioned in the introduction and at the beginning of Section \ref{proofs}, 
   Theorem \ref{thm1} was proved for Riemannian metrics in \cite{japan,archive}. We present an alternative proof, which is much shorter (modulo the results of the previous sections and a  nontrivial result of Tanno  \cite{Tanno}).

  We assume that  $g$   is a complete  Riemannian 
  metric  on a connected manifold such that its degree of mobility is $\ge 3$.   
   Then, by Corollary \ref{ta}, the function  $\lambda$ is a solution of \eqref{tanno}.
   If the metrics are not  affine equivalent, $\lambda$ is not identically constant. 
   
   Let us first assume that the constant  $B$ in the equation \eqref{tanno} is negative. 
   Under this assumption,  
   the equation  \eqref{tanno} was studied by Obata \cite{Obata},  Tanno \cite{Tanno}, and Gallot   \cite{G}. 
    Tanno \cite{Tanno} and Gallot \cite{G} proved  that a complete Riemannian $g$  such that there exists a nonconstant function $\lambda $ satisfying \eqref{tanno} must have a constant positive sectional  curvature.  Applying this  result in our situation, we obtain the claim.

   Now, let us suppose $B\ge 0$. Then,  one can slightly modify the proof from Section \ref{proof1} to obtain the claim. More precisely, substituting \eqref{un5} in \eqref{phi0}, 
   we obtain the following  analog of the equation \eqref{un3}:  
   \begin{equation} \dddot p = 4  B g(\dot\gamma, \dot\gamma) \dot p \label{un3a}.  
   \end{equation}  
   If $B=0$, the equation coincides with \eqref{un3}. Arguing as in Section \ref{proof1}, we obtain that $\phi$ is constant along the geodesic. 
   
   If $B>  0$,   the   general solution  of the equation \eqref{un3a} is 
 \begin{equation}C + C_+e^{2\sqrt{B g(\dot\gamma, \dot\gamma)}\cdot t}+ C_-e^{-2\sqrt{B g(\dot\gamma, \dot\gamma)}\cdot  t}\label{un10} . \end{equation}
  Then, the function $\tau$ satisfies the ODE $\dot\tau = \tfrac{1}{C + C_+e^{2\sqrt{B g(\dot\gamma, \dot\gamma)}\cdot t}+ C_-e^{-2\sqrt{B g(\dot\gamma, \dot\gamma)}\cdot t}} $  implying 
  
  \begin{equation} \tau(t)  =  \int_{t_0}^t \tfrac{d\xi}{C + C_+e^{2\sqrt{B g(\dot\gamma, \dot\gamma)}\cdot \xi }+ C_-e^{-2\sqrt{B g(\dot\gamma, \dot\gamma)}\cdot \xi}} \ \  + \const . \label{un7} \end{equation}

  If one of the constants $C_+, C_-$ is  not zero,  the integral \eqref{un7}  is  bounded from one side, or explodes in finite time.  In both cases,  $\tau $ is not   a diffeomorphism of   $\mathbb{R}$ on itself, i.e., one  of the metrics is not complete. 
  The only possibility for $\tau $ to be a diffeomorphism of   $\mathbb{R}$ on itself is  $C_+=C_-=0$. Finally, $\phi$ is a constant along the  geodesic $\gamma$.

    Since every two points  of a connected complete Riemannian manifold 
     can be connected  by  a geodesic, $\phi$ is a constant, so that $\phi_i\equiv 0$, 
     and the metrics are affine equivalent by \eqref{c1}. \qed

  \begin{Rem}
  Similar idea (contracting the equation with lightlike geodesic and investigating the obtained ODE along the geodesic) was recently used in \cite{conformal_einstein,cones}
  \end{Rem}

\subsection{Proof of Theorem 2} \label{thelast}

Let $g$ be a complete pseudo-Riemannian metric on a connected  closed  manifold  $M^n$ such  that for no $\const\ne 0$ the metric $\const\cdot g$ is Riemannian (if $g$ is Riemannian, Theorem \ref{thm2} follows from Theorem \ref{thm1}).  We assume that the degree of mobility of $g$ is $\ge 3$. Our 
goal is to show that every   metric $\bar g$ geodesically equivalent to $g$ is actually affine equivalent to $g$.

 We consider the function $\lambda$  constructed by \eqref{lam}  for the  solution $a_{ij}$ of \eqref{basic} given by \eqref{a}.  We consider a lightlike geodesic $\gamma(t)$ of the metric $g$, and the function $\lambda(\gamma(t))$. 
By Corollary \ref{temp3}, 
the function $\lambda(\gamma(t))$ satisfies the ODE $\tfrac{d^3}{dt^3} \lambda(\gamma(t))= 0$.  
Hence $\lambda(\gamma(t))= C_2t^2 + C_1t + C_0$. If $C_2\ne 0$, or $C_1\ne 0$, then the function $\lambda$ is not bounded; that contradicts the compactness of  the manifold. Thus $\lambda(\gamma(t))$ is constant along every  lightlike geodesic. Since every two points can be connected by a sequence of lightlike geodesics, $\lambda $ is constant. Thus $\lambda_i=0$ implying in view of  (\ref{lambda}) that $\phi_i=0$ implying in view of \eqref{LC}  that the metrics are affine equivalent.   \qed

\end{document}